\documentclass[12pt]{amsart}

\usepackage{graphicx}
\usepackage{amsfonts}
\usepackage{epsf}
\usepackage{amssymb}
\usepackage{amsmath}
\usepackage{amscd}

\newtheorem{theorem}{Theorem}[section]
\newtheorem{proposition}[theorem]{Proposition}
\newtheorem{lemma}[theorem]{Lemma}
\newtheorem{definition}[theorem]{Definition}

\newtheorem{corollary}[theorem]{Corollary}
 
\newtheorem{remark}[theorem]{Remark} 
\newtheorem{example}[theorem]{Example}

\newcommand{\kerr}{\mbox{Ker} }

\newcommand{\F}{\mathbb{F}}
\newcommand{\lr}{\langle \cdot , \cdot \rangle}

\input xy
\xyoption{all}

\begin{document}

\title{Rational Witt classes of pretzel knots}
\begin{abstract}
In his two pioneering articles \cite{levine1, levine2} Jerry Levine introduced and completely determined the 
algebraic concordance groups of odd dimensional knots. He did so by defining a host of invariants of algebraic concordance which he showed were a complete set of invariants. While being very powerful, these  invariants are in practice often hard to determine, especially for knots with Alexander polynomials of high degree. We thus propose the study of a weaker set of invariants of algebraic concordance -- the rational Witt classes of knots. Though these are rather weaker invariants than those defined by Levine, they have the advantage of lending themselves to quite manageable computability. We illustrate this point by computing the rational Witt classes of all pretzel knots. 
We give many examples and provide applications to obstructing sliceness for pretzel knots. Also, we obtain explicit formulae 
for the determinants and signatures of all pretzel knots.  

This article is dedicated to Jerry Levine and his lasting mathematical legacy; on the occasion of the conference 
\lq\lq Fifty years since Milnor and Fox\rq\rq \, held at Brandeis University on June 2--5, 2008. 
\end{abstract}
\author{Stanislav Jabuka}
\address{Department of Mathematics and Statistics,  
University of Nevada, Reno, NV 89557}
\email{jabuka@unr.edu}
\maketitle
\section{Introduction}
%
%
%
\subsection{Preliminaries}
In his seminal papers \cite{levine1, levine2} Jerry Levine introduced and determined the algebraic concordance groups $\mathcal C _n$ of concordance classes of embeddings of $S^n$ into $S^{n+2}$. These groups had previously been found by Kervaire \cite{kervaire} to be trivial for $n$ even; for $n$ odd, Levine proved that\footnote{For brevity, we denote the infinite direct sum $\oplus _{i=1}^\infty \mathbb{Z}_p$ simply by $\mathbb{Z}_p^\infty$ hoping the reader will not confuse the latter with the product of an infinite number of copies of $\mathbb{Z}_p$. Throughout the article, $\mathbb{Z}_p$ denotes $\mathbb{Z}/p\mathbb{Z}$.}
$$\mathcal C_n \cong   \mathbb{Z}^\infty \oplus  \mathbb{Z}_2^\infty  \oplus  \mathbb{Z}_4^\infty $$ 
Levine achieved this remarkable result by considering a natural homomorphism $\varphi_n : \mathcal C_n \to \mathcal I(\mathbb{Q})$ from the algebraic concordance group $\mathcal C_n$ into the concordance group of isometric structures $\mathcal I(\mathbb{Q})$ on finite dimensional vector spaces over $\mathbb{Q}$ (we describe
$\mathcal I(\mathbb{Q})$ in detail in section \ref{isostr} below). He 
constructed a complete set of invariants of concordance of isometric structures and used these invariants to 
show that $\mathcal I(\mathbb{Q}) \cong  \mathbb{Z}^\infty \oplus  \mathbb{Z}_2^\infty  \oplus  \mathbb{Z}_4^\infty $. Moreover, he showed that he map $\varphi_n: \mathcal C_n \to \mathcal I(\mathbb Q)$ is injective and that 
its image is large enough to itself contain a copy of 
$\mathbb{Z}^\infty \oplus  \mathbb{Z}_2^\infty  \oplus  \mathbb{Z}_4^\infty$, thereby establishing the isomorphism 
$\mathcal C_n \cong \mathbb{Z}^\infty \oplus  \mathbb{Z}_2^\infty  \oplus  \mathbb{Z}_4^\infty$. In this article we
focus exclusively on the case of $n=1$. 
 
To determine the values of Levine's complete set of invariants for a given knot $K$, one is required to find the irreducible symmetric factors of the Alexander polynomial $\Delta _K(t)$ of $K$. As the 
question of whether or not a given polynomial is irreducible is a difficult one in general, the task of determining 
all the irreducible factors of a given symmetric polynomial can be quite intractable, more so as the degree of the polynomial grows. To circumnavigate this issue, 
we consider another homomorphism $\varphi :\mathcal C_1 \to W(\mathbb{Q})$ from the algebraic concordance group $\mathcal C_1$  into the Witt ring over the rationals ($W(\mathbb Q)$ is described in detail in section \ref{wittrings}, for a 
brief description see section \ref{results} below). 
The isomorphism type of $W(\mathbb{Q})$ as an Abelian group is well understood and is given by  
$W(\mathbb{Q}) \cong \mathbb{Z} \oplus \mathbb{Z}_2^\infty \oplus \mathbb{Z}_4^\infty$. 
The maps $\varphi$ and $\varphi _1$ fit into the commutative diagram 

\centerline{
\xymatrix{
\mathcal C _1 \ar[rr]^{\varphi_1} \ar[dr]_{\varphi }  &   & \mathcal I(\mathbb{Q})  \ar[dl]^{\psi } \\
 &  W(\mathbb{Q})  &   \\
}
}
\vskip1mm
\noindent 
From simply knowing the isomorphism types of $\mathcal C_1$ and $W(\mathbb{Q})$, it is clear 
that $\varphi :\mathcal C_1 \to W(\mathbb{Q})$ cannot be injective and a loss of information must occurs in passing from $K\in \mathcal C_1$ to $\varphi(K) \in W(\mathbb{Q})$. The payoff being that one is no longer required to factor polynomials. Indeed, to determine $\varphi (K)$ for a given knot $K\subset S^3$  
one only needs to use the Gram-Schmidt orthogonalization process along with a simple \lq\lq reduction\rq\rq argument (described in section \ref{GSRed}). The Gram-Schmidt process is completely algorithmic (in contrast with polynomial factorization) and is readily available in many mathematics software packages.  
\vskip3mm
To goal of this article then is to underscore the computability and usefulness of the rational Witt classes $\varphi(K)$. Their determination is almost entirely algorithmic and often straightforward, if tedious, to calculate. We illustrate our point by focusing on a concrete family of knots --  the set of pretzel knots. This family is large enough to reflect a number of varied properties of the invariant $\varphi$ and yet tractable enough so that a complete determination of the rational Witt classes is possible. We proceed by giving a few details about pretzel knots first and then state our main results. 
\subsection{Statement of results} \label{results}
\begin{figure}[htb!] 
\centering
\includegraphics[width=8cm]{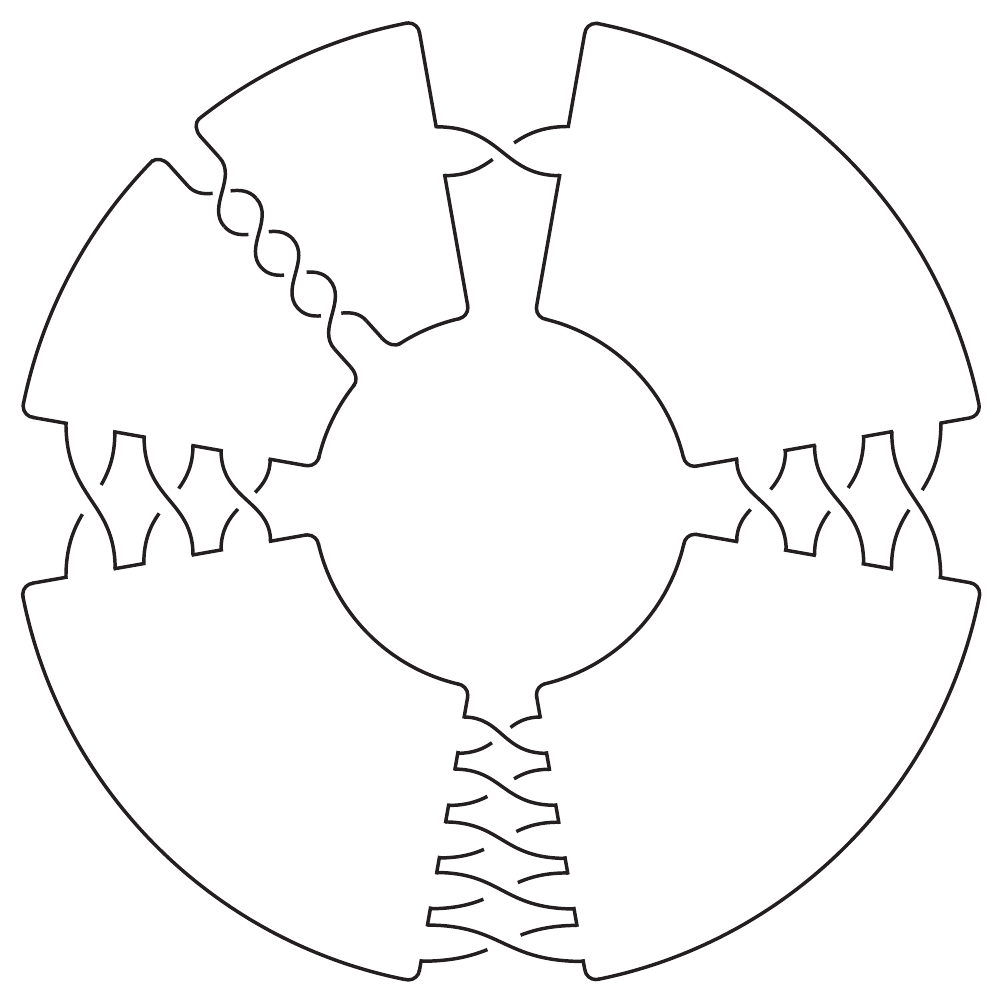}
\caption{The pretzel knot $P(-1,3,-5,3,4)$. }  \label{pic1}
\end{figure}

Given a positive integer $n$ and integers $p_1,p_2,...,p_n$,  let $P(p_1,p_2,...,p_n)$ denote the $n$-stranded pretzel knot/link. It is obtained by taking $n$ pairs of parallel strands, introducing $p_i$  half-twists into the $i$-th strand and capping the strands off by $n$ pairs of bridges. The signs of the $p_i$ determine the handedness of the corresponding half-twists. Our convention is that $p_i>0$ corresponds to right-handed half-twists, see figure \ref{pic1} for an example. We limit our considerations to knots and moreover require that $n\ge 3$ and that $p_i\ne 0$ (the purpose of these two limitations is to exclude connected sums of torus knots/links). There are 3 categories of choices of the parameters $n,p_1,...p_n$ which lead to knots, namely
\begin{align} \label{categories}
\mathbf{(i)} & \mbox{ $n$ is odd and all exept one of the $p_i$ are odd.  } \cr
\mathbf{(ii)} & \mbox{ $n$ is even and all exept one of the $p_i$ are odd.  } \cr
\mathbf{(iii)} & \mbox{ $n$ is odd and all $p_i$ are odd.  } 
\end{align}
As we shall see, these categories exhibit slightly different behavior as far as their images in $W(\mathbb{Q})$. 
Pretzel knots are invariant under the action of $\mathbb{Z}_n$ by cyclic permutation, i.e. 
$P(p_1,p_2,...,p_{n-1},p_n) = P(p_n,p_1,p_2,,...,p_{n-1})$. 
We use this symmetry to fix the convention that if $P(p_1,...,p_n)$ comes from either category (i) or (ii) 
above, we let $p_n$ be the unique even integer among $p_1,...,p_n$. 

To state our results, we need to give a brief description of the rational Witt ring $W(\mathbb{Q})$, a more copious 
exposition is provided in section \ref{wittrings}. As a set, $W(\mathbb{Q})$ consists of equivalence classes of pairs 
$(\langle \cdot , \cdot \rangle , V)$ 
where $V$ is a finite dimensional vector space over $\mathbb{Q}$ and $\langle \cdot , \cdot \rangle :V\times V\to 
\mathbb{Q}$ is a non-degenerate symmetric bilinear form. We say that a pair 
$(\langle \cdot , \cdot \rangle , V)$ is {\em metabolic} or {\em totally isotropic} is there exits a half-dimensional 
subspace $W\subset V$ such that $\langle \cdot, \cdot \rangle |_{W\times W} \equiv 0$. We will be adding 
pairs $(\langle \cdot , \cdot \rangle_1 , V_1)$ and $(\langle \cdot , \cdot \rangle_2 , V_2)$ by direct summing them, thus
$$(\langle \cdot , \cdot \rangle_1 , V_1) \oplus (\langle \cdot , \cdot \rangle_2 , V_2) = (\langle \cdot , \cdot \rangle_1\oplus \langle \cdot , \cdot \rangle_2 , V_1\oplus V_2)$$
With this understood, the equivalence relation on $W(\mathbb{Q})$ is the one by which 
$(\langle \cdot , \cdot \rangle_1 , V_1)$ is equivalent to $(\langle \cdot , \cdot \rangle_2 , V_2)$ 
if $(\langle \cdot , \cdot \rangle_1 , V_1)  \oplus ( - \langle \cdot , \cdot \rangle_2 , V_2)$ is metabolic. One proceeds 
to check that addition is commutative and indeed well defined on $W(\mathbb{Q})$, giving $W(\mathbb{Q})$ the
structure of an Abelian group. 

It is not hard to obtain an explicit presentation of
$W(\mathbb{Q})$ (see theorem \ref{wittpresentation} in section \ref{wittrings}), for now however it will suffice to point out that $W(\mathbb{Q})$ is generated by the set  $\{ \langle a \rangle   \in W(\mathbb{Q}) \, | \, a\in \mathbb{Q}-\{0\} \, \}$. 
Here $\langle a \rangle$ stands for $(\lr_a,\mathbb{Q})$ where $\lr_a$ is the form on $\mathbb{Q}$ specified by 
$\langle 1, 1 \rangle _a=a $. 
\vskip2mm

Given a knot $K\subset S^3$, pick an oriented, genus $g$ Seifert surface $\Sigma \subset S^3$ and consider the linking pairing 
$\ell k : H_1(\Sigma ;\mathbb{Z}) \times H_1(\Sigma  ; \mathbb{Z}) \to \mathbb{Z}$ given by 
$$ \ell k (\alpha, \beta) = \mbox{ linking number between $\alpha$ and $\beta ^+$ } $$
where $\beta^+$ is a small push-off of $\beta$ in the preferred normal direction of $\Sigma $ determined 
by its orientation. Extending $\ell k$ to $H_1(\Sigma ;\mathbb{Q})$ linearly and letting 
$\langle \cdot , \cdot \rangle : H_1(\Sigma ;\mathbb{Q})\times H_1(\Sigma  ; \mathbb{Q}) \to \mathbb{Q}$ 
be $\langle \alpha, \beta \rangle = \ell k ( \alpha,\beta) + \ell k (\beta, \alpha)$, defines a non-degenerate symmetric bilinear pairing on the 
rational vector space $H_1(\Sigma_g;\mathbb{Q})$. We use this to define 
$$\varphi (K) = ( \langle \cdot , \cdot \rangle, H_1(\Sigma ;\mathbb{Q})) \in W(\mathbb{Q})$$
which we refer to as the {\em rational Witt class of $K$}. 
According to \cite{levine1}, $\varphi(K)$ is well defined and only depends on $K$ (as an oriented knot) but not 
on the particular choice of Seifert surface $\Sigma $. In fact, $\varphi (K)$ only depends on the algebraic 
concordance class of $K$. 
\vskip3mm
With these descriptions out of the way, we are now ready to state our main results. 
%
%
%
\begin{theorem}  \label{main1}
Consider category $(i)$ from \eqref{categories}, i.e. let $n\ge 3$ be an odd integer, let $p_1,...,p_{n-1}$ be odd integers and let $p_n\ne 0$ be an even integer. 
Then the rational Witt class of the pretzel knot $P(p_1,...,p_n)$ is given by
\begin{align} \nonumber
\varphi(P(p_1,...,p_n))  & = \bigoplus _{i=1}^{n-1}  \left(\left\langle s_i \cdot  1\cdot 2\right\rangle  \oplus  
\left\langle s_i  \cdot 2\cdot 3\right\rangle \oplus  ...\oplus  \left\langle s_i \cdot (|p_i|-1)\cdot |p_i|) \right\rangle \right) \oplus \cr 
&\quad \quad \quad \quad   \oplus \langle -(p_1\cdot ... \cdot p_{n-1})\cdot \det P(p_1,...,p_{n-1}) \rangle \oplus \cr 
& \quad \quad \quad \quad   \oplus \langle \det P(p_1,...,p_{n-1}) \cdot \det P(p_1,...,p_n) \rangle  
\end{align}
where $s_i = -Sign(p_i)$. 
%
%
%
%
The two determinants appearing above equal 
\begin{align} \nonumber
\det P(p_1,...,p_n) & = \prod _{i=1}^{n} p_1 \cdot ... \cdot \hat p_i \cdot 
... \cdot p_{n} \cr 
\det P(p_1,...,p_{n-1}) & = \prod _{i=1}^{n-1} p_1 \cdot ... \cdot \hat p_i \cdot ... \cdot p_{n-1}
\end{align} 
\end{theorem}
As is customary in the literature, having a hat decorate a variable in a product indicates that the factor 
should be left out. For example $p_1 \cdot \hat p_2 \cdot p_3$ stands for $p_1 \cdot p_3$.
\begin{theorem} \label{main2}
Consider category $(ii)$ from \eqref{categories}, that is, let $n\ge 3$ be an even integer, let $p_1,...,p_{n-1}$ be odd integers and let $p_n\ne 0$ be an even integer. 
Then  the rational Witt class of the pretzel knot $P(p_1,...,p_n)$ is  
\begin{align} \nonumber
\varphi(P(p_1,...,p_n))  & = \bigoplus _{i=1}^{n}  \left(\left\langle s_i \cdot  1\cdot 2\right\rangle  \oplus  
\left\langle s_i  \cdot 2\cdot 3\right\rangle \oplus  ...\oplus  \left\langle s_i \cdot (|p_i|-1)\cdot |p_i|) \right\rangle \right) \oplus \cr 
&\quad \quad \quad   \oplus \langle - (p_1\cdot ... \cdot p_n) \cdot  \det P(p_1,...,p_{n}) \rangle
\end{align}
where $s_i = -Sign(p_i)$ and the determinant $\det P(p_1,...,p_n)$ is again given by 
$$\det P(p_1,...,p_n) = \prod _{i=1}^n p_1 \cdot ... \cdot \hat p_i \cdot ... \cdot p_n$$ 
\end{theorem}
To state the next theorem we introduce some auxiliary notation first: Let $\sigma _j(t_1,...,t_m)$ denote the degree $j$ 
symmetric polynomial in the variables $t_1, ..., t_m$. For example, $\sigma _1(t_1,...,t_m) = t_1+...+t_m$ 
while $\sigma _m(t_1,...,t_m) = t_1 \cdot ... \cdot t_m$. We adopt the convention that $\sigma _0 (t_1,...,t_m) = 1$. 
With this in mind, we have
\begin{theorem} \label{main3}
Consider category $(iii)$ from \eqref{categories}. Thus, let $n\ge 3$ and $p_1,...,p_n$ be odd integers and 
let $\sigma _i$ stand as an abbreviation for the integer $\sigma _i(p_1,...,p_{i+1})$. 
Then  the rational Witt class of the pretzel knot $P(p_1,...,p_n)$ is given by 
\begin{align} \nonumber
\varphi(P(p_1,...,p_n))  & =\langle \sigma _0 \cdot \sigma _1 \rangle  \oplus  
\langle \sigma _1 \cdot \sigma _2\rangle \oplus  ...\oplus  \langle \sigma _{n-2} \cdot \sigma _{n-1} \rangle
\end{align}
We note that $\det P(p_1,...,p_n) = \sigma _{n-1}$. 
\end{theorem}
\begin{remark}
To put the results of theorems \ref{main1} -- \ref{main3} into perspective, we would like to point out that at the time of this writing, the algebraic concordance orders aren't known yet even for 
the 3-stranded pretzel knots $P(p_1,p_2,p_3)$ from category $(i)$ in \eqref{categories}. The chief reason for this is that this family contains knots with Alexander polynomials of arbitrarily high degree. 

In contrast, the algebraic concordance orders of $P(p_1,p_2,p_3)$ coming from category $(iii)$ in \eqref{categories} are well understood and follow easily from Levine's article \cite{levine2}, see remark \ref{levineremark} below. All non-trivial knots in this family are of Seifert genus $1$. 
\end{remark}
%
%
%
\subsection{Applications and examples} \label{appandex}
While theorems \ref{main1} -- \ref{main3} give $\varphi(K)$ in terms of the generators of $W(\mathbb{Q})$, in 
concrete cases one can determine $\varphi (K)$ as a specific element in $\mathbb{Z}\oplus \mathbb{Z}_2^\infty \oplus \mathbb{Z}_4^\infty \cong W(\mathbb{Q})$. We give a host of examples of this nature next. 
Such computations rely on an understanding of the isomorphism between $W(\mathbb{Q})$ and $\mathbb{Z}\oplus \mathbb{Z}_2^\infty \oplus \mathbb{Z}_4^\infty$. This isomorphism is completely explicit and easily computed, we explain it in some detail in section \ref{wittrings}. For now we merely present the results of our computations, the full details are deferred to section \ref{examples-section}.

After presenting a several concrete examples, we turn to general type corollaries of theorems \ref{main1} -- \ref{main3}. The ultimate goal of course is to have a set of numerical conditions on $n,p_1,...,p_n$ which would pinpoint
the order of $\varphi(K)$ in $W(\mathbb{Q})$. The obstacle to achieving this is number theoretic in 
nature and we have been unable to overcome it in its full generality. However, we are able to give such conditions for the case of $n=3$ and for some special cases when $n\ge 4$.   
\vskip3mm
As we shall see in section \ref{wittrings}, a necessary condition for $\varphi(K)$ to be zero in $W(\mathbb{Q})$ is
that $\sigma(K)=0$ and $|\det K| = m^2$ for some odd integer $m$. If only the first of these conditions holds, then 
$\varphi(K)$ is at least of order $2$ in $W(\mathbb{Q})$. With this in mind the next examples testify that the 
rational Witt classes carry significantly more information than merely the signature and determinant. 
We start with a useful definition
\begin{definition} \label{upstab}
If $p$ is an odd integer, we shall say that the knot 
$$P(p_1,...,p_{i-1},p,p_i,...,p_{j-1},-p,p_j,...,p_n)$$
is gotten from $P(p_1,...,p_n)$ by an upward stabilization (or conversely that $P(p_1,...,p_n)$ is obtained 
from $P(p_1,...,p_{i-1},p,p_i,...,p_{j-1},-p,p_j,...,p_n)$ by a downward stabilization). 
\end{definition}
\begin{example} \label{example1}
Let $K_1$, $K_2$ and $K_3$ be the knots 
$$K_1=P(21,13,-17,-15,12) \quad \quad K_2=P(-3,-3,-7,5,2) \quad \quad K_3=P(-3,-5,7,9,6)$$
from category $(i)$ and let $K=K_1\#K_2\#K_3$. The $\sigma (K)=0$ but $\varphi (K)$ has order $4$ 
in $W(\mathbb{Q})$. Thus $K$ has concordance order at least $4$. The same holds if 
$K_i$ is replaced by a knot gotten from $K_i$ by any finite number of upward stabilizations. 
\end{example}
\begin{example} \label{example2}
Let $K_1$ and $K_2$ be the knots 
$$K_1=P(7,3,-5,2) \quad \quad K_2=P(-19,-15,21,10)$$
from category $(ii)$ and let $K=K_1\#K_2$. The $\sigma (K)=0$ but $\varphi (K)$ has order $4$ 
in $W(\mathbb{Q})$ and therefore also in concordance group. The same is true if $K_i$ is replaced by a 
knot gotten from $K_i$ by any finite number of upward stabilizations.   
\end{example}
\begin{example} \label{example3}
Let $K$ be a knot obtained by a finite number of upward stabilization from either 
$$P(-3,9,15,-5-5) \quad \mbox{ or } \quad P(-3,-5,-11,15,15)$$
from category $(iii)$. Then the signature of $K$ is zero, the determinant of $K$ is a square but $\varphi(K) \ne 0 \in W(\mathbb{Q})$. Consequently, no such $K$ is slice. 
\end{example}
\begin{example} \label{example4}
Let $K_1$, $K_2$ and $K_3$ be the knots  
$$K_1=P(21,13,-17,-15,12) \quad \quad K_2=P(-19,-15,21,10) \quad \quad K_3=P(-15,-7,-7,13,11)$$
from the categories $(i)$, $(ii)$ and $(iii)$ and let $K=K_1\#K_2\#K_3$. Then $\sigma (K) = 0$ 
but $\varphi (K)$ is of order $4$ in $W(\mathbb{Q})$. The same holds under replacement of $K_i$ by 
upward stabilizations. 
\end{example}
The details of the above computations can be found in section \ref{examples-section}. 
We now turn to more general corollaries of theorems \ref{main1} -- \ref{main3}. 
\begin{theorem} \label{coro1}
Consider a 3-stranded pretzel knot $K=P(p,q,r)$ with $p,q,r$ odd. Then the order of $\varphi (K)$ in $W(\mathbb{Q})$ 
is as follows:
\begin{itemize}
\item $\varphi(K)$ is or order 1 in $W(\mathbb{Q})$ if and only if $\det K=-m^2$ for some odd $m\in \mathbb{Z}$.
\item $\varphi(K)$ is of order 2 in $W(\mathbb{Q})$ if and only if $\det K<0$, $\det K$ is not a square and $|\det K| \equiv 1\, \,  (\mbox{mod } 4)$.
\item $\varphi(K)$ is of order 4 in $W(\mathbb{Q})$ if and only if $\det K<0$ and $|\det K| \equiv 3\, \,  (\mbox{mod } 4)$.
\item $\varphi(K)$ is of infinite order $W(\mathbb{Q})$ if and only if $\det K >0$.
\end{itemize}
Recall that $\det K = pq+pr+qr$. 
\end{theorem}
\begin{theorem}  \label{oddoddeven-theorem}
Consider again $K=P(p,q,r)$ but with $p,q$ odd and with $r\ne 0$ even. Then 
$\varphi(K)$ is of finite order in $W(\mathbb{Q})$ if and only if 
$$p+q=0\quad \quad \quad \mbox{ or } \quad \quad \quad p+q=\pm2 \, \, \mbox{ and } \, \, \det K >0 $$  
The order of $\varphi(K)$ in $W(\mathbb{Q})$ in these cases is as follows:
\begin{itemize}
\item If $p+q=0$ then $\varphi(K)$ has order $1$ in $W(\mathbb{Q})$. 
\item If $p+q=\pm 2$ and $\det K>0$ then 
\begin{itemize}
\item  $\varphi(K)$ is of order $1$ in $W(\mathbb{Q})$ if $\det K=m^2$ for some odd integer $m$.
\item  $\varphi(K)$ is of order $2$ in $W(\mathbb{Q})$ if $\det K$ is not a square and is congruent to
$1 \, (\mbox{mod } 4)$. 
\item  $\varphi(K)$ is of order $4$ in $W(\mathbb{Q})$ if $ \det K \equiv 3 \, (\mbox{mod } 4)$. 
\end{itemize} 
\end{itemize}
Here too $\det K = pq+pr+qr$.
\end{theorem}
A slightly more general version of this theorem is given in theorem \ref{oddoddeven-theorem-general}.  
\begin{remark} \label{levineremark}
As already mentioned above, the algebraic concordance orders of the knots $P(p,q,r)$ with $p,q,r$ odd 
are known by work of Levine \cite{levine2} and agree with the orders of $\varphi (P(p,q,r))$ in $W(\mathbb{Q})$. 
The analogues of the results of theorem \ref{oddoddeven-theorem} are not known for the algebraic concordance group. However, according to theorem \ref{main5} below, it is clear that when $r$ is even, the order of $\varphi(P(p,q,r))$ in $W(\mathbb{Q})$ and the order of $P(p,q,r)$ in $\mathcal C_1$ are different in general.
\end{remark}

\begin{remark}
The condition on the congruency class mod $4$, appearing in both theorems \ref{oddoddeven-theorem} and \ref{coro1}, is reminiscent of a similar condition appearing in a beautiful (and much stronger) theorem by 
Livingston and Naik \cite{chuckswatee}: If $K$ is a knot with $\det K = \wp\cdot \beta$ where $\wp$ is a prime
congruent to $3$ mod $4$ and $\gcd (\wp, \beta) = 1$, then $K$ has infinite order in the topological concordance group.    
\end{remark}
\begin{theorem} \label{coro2}
Consider a pretzel knot $K=P(p_1,...,p_n)$ from category $(i)$ in \eqref{categories}, i.e. assume that $n,p_1,....,p_{n-1}$ are odd, $n\ge 3$ and $p_n\ne 0$ is even.  
Additionally, suppose that the $p_1,...,p_{n-1}$ are all mutually coprime. Then $\varphi (K) = 0\in W(\mathbb{Q})$ if and only if $\sigma(K)=0$ and $\det K = \pm m^2$ for some odd  $m\in \mathbb{Z}$. 
\end{theorem}
Seeing as the torsion subgroups of $\mathcal C_1$ and $W(\mathbb{Q})$ are isomorphic, one can't help but 
speculate whether $\varphi |_{Tor (\mathcal C_1)}:Tor (\mathcal C_1) \to W(\mathbb{Q})$ is injective. Unfortunately this is not the case as the next theorem testifies. 
\begin{theorem} \label{main5}
Consider the knot $K=P(5,-3,8)$. All Tristram-Levine signatures $\sigma _\omega (K)$ vanish but 
$K$ is not trivial in $\mathcal C_1$. On the other hand, the rational Witt class $\varphi (K)$ is zero . 
Thus, $K$ is a nontrivial element of $\kerr (\varphi)\cap Tor(\mathcal C_1)$. 
\end{theorem}
\begin{remark} \label{lowcrossings}
We would like to point out that for knots $K$ with 10 or fewer crossings, $K$ is algebraically slice if and only if 
$\varphi(K)$ is zero in $W(\mathbb{Q})$. This follows by inspection, using KnotInfo\footnote{A web site 
created by Chuck Livingston and maintained by Chuck Livingston and Jae Choon Cha. The site contains a wealth of information about knots with low crossing number. It can be found at {\tt http://www.indiana.edu/$\sim~$\hspace{-1mm}knotinfo}.}, and relying on the fact that if $\varphi(K)=0$ then $\sigma (K) = 0$ and $\det K = \pm m^2$.
\end{remark}
As a byproduct of our computations we obtain closed formulae for the signature and determinants 
of all pretzel knots. The formulae for the determinants have already been stated in theorems \ref{main1} -- \ref{main3}, the signature formulae are the content of the next theorem. While these are not directly relevant to our discussion, we list them here in the hopes that they may be useful elsewhere. 
\begin{theorem} \label{signdet}
Let $K=P(p_1,...,p_n)$ be a pretzel knot from either of the 3 categories $(i)$--$(iii)$ from \eqref{categories}. As usual, we assume that $n\ge3$. Then the signature $\sigma(K)$ of $K$ can be computed as follows: 
\begin{enumerate}
\item If $n,p_1,...,p_{n-1}$ are odd and $p_n\ne 0$ is even, then 
\begin{align} \nonumber
\sigma (K) &= -\left( \sum_{i=1}^{n-1} Sign(p_i)\cdot (|p_i|-1) \right) + Sign(p_1\cdot ... \cdot p_{n-1}\cdot \det P(p_1,...,p_{n-1})) + \cr
& \quad \quad \quad \quad \quad \quad \quad \quad\quad \quad+  Sign(\det P(p_1,...,p_{n-1}) \cdot \det P(p_1,...,p_n)) 
\end{align}
The determinants $\det P(p_1,...,p_n)$ and $\det P(p_1,...,p_{n-1})$ are computed as in theorem \ref{main1}.
\item If $n,p_n$ are even, $p_1,...,p_{n-1}$ are odd and $p_n\ne 0$, then 
\begin{align} \nonumber
\sigma (K) &= - \left( \sum_{i=1}^{n-1} Sign(p_i)\cdot (|p_i|-1) \right) + Sign(p_1\cdot ... \cdot p_{n} \cdot \det P(p_1,...,p_{n})) 
\end{align}
where $\det P(p_1,...,p_n)$ is as computed in theorem \ref{main2}.
\item If $n,p_1,...,p_n$ are all odd, then 
\begin{align} \nonumber
\sigma (K) &= \sum _{i=1}^{n-1} Sign(\sigma_{i-1} \cdot \sigma _i) 
\end{align}
where $\sigma_i = \sigma _i(p_1,...,p_{i+1})$ as in theorem \ref{main3}. 
\end{enumerate}
\end{theorem}
For example, if $K=P(p_1,...,p_n)$ with $n,p_1,...,p_n$ odd and $p_i>0$ for all $i$, then $\sigma _i>0$ 
for all $i$ also and therefore $\sigma (K) = n-1$. As another example consider the case of $n,p_n$ even 
and $p_1,...,p_{n-1}$ odd and again $p_i>0$ for all $i$. Then $\sigma (K) = n+1-(p_1+...+p_n)$. 
\subsection{Organization} 
Section \ref{acg} provides background on the three flavors of algebraic concordance groups $\mathcal C_1$, $\mathcal I(\mathbb{Q})$ and $W(\mathbb{Q})$ encountered in the
introduction. The relationships between these groups are also made more transparent. In section \ref{linking-section} the first steps towards computing $\varphi (P(p_1,...,p_n))$ are taken in that  specific Seifert surfaces are picked for the 
knots along with specific bases for their first homology. These choices allow us to determined a linking 
matrix for the knots. Section \ref{diagonalizing-section} explains how one can diagonalize the linking matrices found in section \ref{linking-section}, leading 
to proofs of theorems \ref{main1}, \ref{main2} and \ref{main3}. More detailed versions of these theorems are 
provided in theorems \ref{diagmat1}, \ref{diagmat2} and \ref{diagmat3} respectively.  Section \ref{examples-section}
is devoted to computations of examples and shows how theorems \ref{main1} -- \ref{main3} imply the results from examples \ref{example1} -- \ref{example4} stated above. The final section provides proofs for theorems
\ref{coro1}, \ref{oddoddeven-theorem}, \ref{coro2} and \ref{main5}. 

{\bf Acknowledgement } In the preparation of this work I have greatly benefitted from conversations 
with Chuck Livingston. I am grateful for his generousity in sharing his insight and expertise. 
\section{Algebraic concordance groups} \label{acg}
In this section we describe the three algebraic concordance groups mentioned in the introduction, namely 
\vskip2mm
\begin{tabular}{rcl}
$\mathcal C_1$ & -- & The algebraic concordance group of classical knots in $S^3$.  \cr
$\mathcal{I} (\mathbb{F})$ & -- & The concordance group of isometric structures over the field $\mathbb{F}$. 
 \cr
$W(\mathbb{\mathbb F})$ & -- & The Witt ring of non-degenerate, symmetric, bilinear forms over $\mathbb{F}$. 
\end{tabular}
\vskip2mm
\subsection{The algebraic concordance group $\mathcal C_1$} \label{acg-section}
This section largely follows the exposition from \cite{levine1} with a slight bias towards a coordinate free 
description. 

Our explanation of the algebraic concordance group $\mathcal C_1$ runs largely in parallel to the
description of the Witt ring $W(\mathbb{Q})$ from the introduction. Thus, we shall consider 
pairs $(\lr,L)$ where $L$ is a finitely generated free Abelian group and $\lr : L\times L \to \mathbb{Z}$ is 
a bilinear pairing with the property that $\lr - \lr ^\tau$ is unimodular. Following Levine \cite{levine2}, we shall 
call such pairs {\em admissible pairs}. Here $\lr ^\tau$ denotes the bilinear 
form 
$$ \langle x,y\rangle ^\tau = \langle y,x\rangle $$ 
Note that $\lr$ is not required to be symmetric nor non-degenerate. 
We will say that $(\lr,L)$ is {\em metabolic} or {\em totally isotropic} if there exists a splitting $L\cong L_1 \oplus L_2$ with $rk \, L = 2 \, (rk \, L_1)$ and  $\lr|_{L_1\times L_1} \equiv 0$. We shall add pairs $(\lr _1, L_1)$ and $(\lr_2,L_2)$ 
by direct summing them, i.e. 
$$(\lr _1, L_1)\oplus (\lr_2,L_2) = (\lr_1\oplus \lr_2,L_1\oplus L_2)$$
With these definitions understood, we define the {\em algebraic concordance group $\mathcal C_1$} to be 
the set of pairs $(\lr, L)$ as above, up to the equivalence relation $\sim$ by which 
$$(\lr _1, L_1) \sim (\lr_2,L_2) \quad \mbox{ if and only if } \quad (\lr _1, L_1)\oplus (-\lr_2,L_2) \mbox{ is metabolic.}$$
We shall refer to this equivalence relation as that of 
{\em algebraic concordance}. Under the operation of direct summing, $\mathcal C_1$ 
becomes an Abelian group. An easy check reveals that the inverse of $(\lr,L)$ is $(-\lr,L)$.
The group $\mathcal C_1$ was introduced by Jerry Levine in \cite{levine1} and its isomorphism type was completely determined by him in \cite{levine2}.  
\vskip2mm
The relation of $\mathcal C_1$ to knot theory is as follows: Let $K$ be a knot in $S^3$ and let $\Sigma \subset S^3$ be an oriented genus $g$ Seifert surface for $K$. We shall view the orientation on $\Sigma $ as being 
given by an normal unit vector field $\vec n$ on $\Sigma$. Recall from the introduction that the linking 
pairing $\ell k : H_1(\Sigma ;\mathbb{Z})\times H_1(\Sigma ;\mathbb{Z}) \to \mathbb{Z}$ is defined by 
$$ \ell k (x,y) = \mbox{ linking number of $x$ and $y^+$ } $$
where, by a customary blurring of viewpoints, we interpret $x$ and $y$ as simple closed curves on $\Sigma$. 
With this in mind, $y^+$ is a small push-off of $y$ in the normal direction of $\Sigma$ determined by $\vec n$. 
It is well known (see e.g. \cite{rolfsen}) that $(\ell k , H_1(\Sigma;\mathbb{Z}))$ is an admissible pair
and therefore the assignment 
$(K,\Sigma)\mapsto (\ell k , H_1(\Sigma;\mathbb{Z})) \in \mathcal C_1$ is well defined. As Levine shows in \cite{levine1}, the algebraic concordance class of $(\ell k , H_1(\Sigma;\mathbb{Z}))$ is independent of $\Sigma$
and by abuse of notation, we shall denote it simply by $K$, hoping that no confusion will arise. Levine 
also shows that if $K_1$ and $K_2$ are (geometrically) concordant as knots then their linking forms are algebraically concordant. This statement applies to both smooth and topological (geometric) concordance.  
\subsection{The Witt ring over the field $\mathbb{F}$} \label{wittrings}
For an excellent introduction to Witt rings we advise the reader to consult \cite{lam}, but see also \cite{milnor}
and \cite{omeara}. The first half of this section is a re-iteration of the description for the Witt ring 
$W(\mathbb{Q})$ over the rational numbers extended to arbitrary fields. 

Let $\mathbb{F}$ be a field and consider pairs $(\lr, V)$ where $V$ is a finite dimensional $\F$-vector 
space and $\lr :V\times V\to \F$ is a symmetric, non-degenerate bilinear pairing. By \lq\lq non-degenerate\rq\rq \, 
we mean that the map $v\mapsto \langle \cdot, v\rangle$ provides an isomorphism from $V$ to $V^\ast$. We call 
a pair $(\lr, V)$ {\em metabolic} or {\em totally isotropic} if there exists a subspace $W\subset V$ with 
$\dim _\F V = 2 \dim _\F W$ and such that $\lr |_{W\times W} \equiv 0$. As in the case of $\F = \mathbb{Q}$, we 
define addition of $(\lr_1, V_1)$ and $(\lr_2, V_2)$ by direct sum
$$ (\lr_1, V_1)  \oplus (\lr_2, V_2) = (\lr_1\oplus \lr _2, V_1\oplus V_2)$$
and we proceed to define the equivalence relation $(\lr_1, V_1) \sim (\lr_2, V_2)$ to mean that 
$(\lr_1, V_1)\oplus (-\lr_2, V_2)$ is metabolic. The set of equivalence classes of pairs $(\lr, V)$  
is denoted by $W(\F)$ and called the {\em Witt ring of $\F$}. It becomes an Abelian group 
under the direct sum operation and a commutative ring with the operation of multiplication given by 
tensor products
$$ (\lr_1, V_1)  \otimes (\lr_2, V_2) = (\lr_1\cdot \lr _2, V_1\otimes _\mathbb{F} V_2)$$
The Witt ring $W(\F)$ was introduced by Witt in \cite{witt} and has found renewed prominence in the theory 
of quadratic forms over fields through the work of Pfister (see for example \cite{pfister1,pfister2}).
\vskip1mm

As is usual in the literature, we will denote $\F - \{0\}$ by $\dot{\F}$. 
Let us recall the notation $\langle a \rangle$ already used in the introduction: Given $a\in \dot{\F}$ we let 
$\langle a \rangle$ denote the non-degenerate symmetric bilinear form $(\langle \cdot , \cdot \rangle_a, \F)$ 
specified by $\langle 1,1 \rangle_a = a$.  Note that 
\begin{equation} \label{cancelsquare}
\langle a \rangle = \langle a \cdot b^2\rangle  \in W(\F) \quad  \forall \, b\in \dot \F 
\quad \quad \mbox{ and } \quad \quad \langle b\rangle \oplus \langle -b\rangle = 0 \in W(\F) \quad \forall b\in \dot \F
\end{equation}
The first of these follows from the fact that $f:(\langle a \rangle , \mathbb{F}) \to (\langle a\cdot b^2\rangle , \mathbb{F})$  given by $f(x) = x\cdot b$ is an isomorphism of forms. The second form is clearly metabolic and thus zero in $W(\F)$. These \lq\lq harmless\rq\rq observations are incredibly useful in computations and we will 
rely on them substantially in our sample calculations in section \ref{examples-section}. 

With this notation in mind, the following theorem can be found in 
\cite{lam}. 
\begin{theorem} \label{wittpresentation}
Let $\langle \cdot , \cdot \rangle$ be a non-degenerate symmetric bilinear form on a finite dimensional 
$\F$-vector space $V$ of dimension $n$. Then there exist scalars $d_1,...,d_n\in \dot \F$ such that 
$$ \langle \cdot , \cdot \rangle = \langle d_1 \rangle \oplus  ... \oplus \langle d_n \rangle \in W(\F)$$
Said differently, $W(\F)$ is generated by the set $\{ \langle a \rangle \, | \, a\in \dot \F\}$.  A presentation 
of $W(\F)$ (as a commutative ring) is obtained from these generators along with the relators
\begin{align} \nonumber
(R1) & \quad \quad \quad \langle 1 \rangle - 1  & \cr
(R2) & \quad \quad \quad \langle a \rangle \cdot  \langle b \rangle - \langle a \cdot b \rangle  & a,b \in \dot \F \cr
(R3) & \quad \quad \quad \langle a \rangle + \langle b \rangle - \langle a+b \rangle \cdot (1+\langle a\cdot b \rangle) & a,b \in \dot \F 
\end{align} 
 In other words, $W(\F)$ is isomorphic to quotient of the free commutative ring generated by  the set 
 $\{ \langle a \rangle \, | \, a\in \dot \F\}$ by the ideal generated by elements of the form as in $(R1)$ -- $(R3)$.
 In $(R1)$, the symbol $1$ denotes the multiplicative unit of $W(\F)$.  
\end{theorem}
 
With this we turn to studying some specific Witt rings. We will chiefly be interested in the cases where $\F$ 
is either $\mathbb Q$ or $\F_\wp$ where the latter will be our notation for the finite field of characteristic $\wp\ge 2$.  The next result can again be found in \cite{lam} and also in \cite{milnor}.

\begin{theorem} 
Let $\wp\in \mathbb{Z}$ be a prime. Then there are isomorphisms of Abelian groups 
$$ W(\F_\wp) \cong \left\{ 
\begin{array}{cl}
\mathbb{Z}_2 & \quad ; \quad \wp = 2 \cr
\mathbb{Z}_2 \oplus \mathbb{Z}_2 & \quad ; \quad \wp \equiv 1 \mbox{ $($mod } 4) \cr
\mathbb{Z}_4 & \quad ; \quad \wp \equiv 3 \mbox{ $($mod } 4)
\end{array}
\right.
$$
The generators of $\mathbb{Z}_2\cong W(\F_2)$ and of $\mathbb{Z}_4\cong W(\F_\wp)$ with 
$\wp \equiv 3 \mbox{ $($mod } 4)$ are given by $\langle 1 \rangle$ while the two copies of $\mathbb{Z}_2$ 
in $W(\F_\wp)$ in the case when $\wp \equiv 1 \mbox{ $($mod } 4)$ are generated by $\langle 1 \rangle $ and 
$\langle a \rangle$ where $a\in \dot \F$ is any non-square element. 
\end{theorem}
The origins of the proof of the next theorem go back to Gauss' work on quadratic reciprocity, it 
was re-discovered by Milnor and Tate \cite{milnor}.
\begin{theorem} 
There is an isomorphism of Abelian groups 
$$ \sigma\oplus \partial : W(\mathbb{Q}) \to \mathbb{Z} \oplus \left(\bigoplus _{\tiny \begin{array}{c} \wp \in \mathbb{N} \cr \wp = \mbox{\tiny prime}\end{array}}  W(\F_\wp) \right) $$
where $\sigma :W(\mathbb{Q})\to \mathbb{Z}$ is the signature function while $\partial : W(\mathbb{Q})\to  \left(\oplus _{\wp}  W(\F_\wp) \right)$ is the direct sum of homomorphisms 
$\partial _\wp :W(\mathbb{Q}) \to W(\F_\wp)$ (with $\wp$ ranging over all primes) described on generators of $W(\F_\wp)$ as follows: Given a rational number $\lambda \ne 0$, write it as $\lambda = \wp ^\ell \cdot \beta $ where  $\ell$ is an integer and $\beta$ a 
rational number whose numerator and denominator are relatively prime to $\wp$. Then 
\begin{equation} \label{delp}
\partial _\wp (\langle \wp ^\ell \cdot \beta \rangle ) = \left\{ 
\begin{array}{cl}
0 & \quad ; \quad \ell \mbox{ is even } \cr
\langle \beta \rangle & \quad ; \quad \ell \mbox{ is odd } 
\end{array}
\right. 
\end{equation}
\end{theorem}
\begin{corollary}
As an Abelian group, $W(\mathbb{Q})$ is isomorphic to $\mathbb{Z}\oplus \mathbb{Z}_2^\infty \oplus \mathbb{Z}_4^\infty$. 
\end{corollary}
\subsection{The concordance group of isometric structures} \label{isostr}
For more details on this section, see \cite{levine2}. 

Let $\F$ be a field, then an {\em isometric structure over $\F$} is a triple 
$(\langle \cdot , \cdot \rangle, T, V)$ consisting of a non-degenerate symmetric bilinear form 
$(\lr , V)$ and a linear operator $T:V\to V$ which is an isometry with respect to $\lr$, i.e. 
$\langle Tv , Tw\rangle = \langle v, w \rangle$ for all $v,w\in V$. 
A triple $(\langle \cdot , \cdot \rangle, T, V)$ shall be called {\em metabolic} or {\em totally isotropic} if there is a
half-dimensioinal $T$-invariant subspace $W\subset V$ for which $\lr |_{W\times W} \equiv 0$.  
Much as in the case of the algebraic concordance group $\mathcal C_1$ and the Witt ring $W(\F)$, isometric structures too are added by direct sum $\oplus$. We define two triples $(\langle \cdot , \cdot \rangle_1, T_1, V_1)$ and $(\langle \cdot , \cdot \rangle_2, T_2, V_2)$ to be equivalent if 
$$ (\langle \cdot , \cdot \rangle_1, T_1, V_1) \oplus (-\langle \cdot , \cdot \rangle_2, -T_2, V_2)$$
is metabolic. With these definitions understood, we define the {\em concordance group of isometric structures $\mathcal I(\F)$} as the set of equivalence classes of triples $(\langle \cdot , \cdot \rangle, T, V)$ as above. 
Not surprisingly, $\mathcal{I}(\F)$ becomes an Abelian group under the operation of direct summing. 
\subsection{Maps between the algebraic concordance groups}
Having defined $\mathcal{C}_1$, $W(\mathbb{\F})$ and $\mathcal I(\mathbb{\F})$, we turn to describing 
some natural maps between them in the case when $\F = \mathbb{Q}$. We start by a lemma proved by Levine in \cite{levine2}.
\begin{lemma} 
Let $(\lr , L)$ be an admissible pair (as in section \ref{acg-section}). Then there exists an admissible pair 
$(\lr ' , L')$ algebraically concordant to $(\lr , L)$ and such that $\lr':L'\times L'\to \mathbb{Z}$ is a 
non-degenerate bilinear form.
\end{lemma}
With this in mind, consider an admissible non-degenerate pair $(\lr, L)$. Given any basis $\mathcal B= \{\alpha _1, ..., \alpha _n\}$ of $L$, let $A$ be the matrix representing $\lr$, that is, set $a_{i,j} = \langle \alpha _i, \alpha _j\rangle$ and let $A=[a_{i,j}]$. We define the maps
$\varphi :\mathcal C_1 \to W(\mathbb{Q})$, $\varphi _1: \mathcal C_1 \to \mathcal I(\mathbb Q)$ and 
$\psi : \mathcal I(\mathbb{Q}) \to W(\mathbb Q)$ as in \cite{levine2}
\begin{align} \nonumber
\varphi(\lr, L) & = ( \lr + \lr^\tau , L\otimes _{\mathbb{Z}} \mathbb{Q} )   \cr
\varphi_1 (\lr, L) &= ( A+A^\tau, -A^{-1}A^\tau, L\otimes _{\mathbb{Z}} \mathbb{Q}) \cr
\psi (\langle \cdot , \cdot \rangle, T, V) & =  (\langle \cdot , \cdot \rangle,  V)
\end{align}
It is not hard to verify that the definition of $\varphi_1$ is independent of the choice of the basis $\mathcal B$ of $L$. It is also easy to verify that, with respect to $\mathcal B$, the matrix $-A^{-1}A^\tau$ defines an isometry on 
$L\otimes _{\mathbb{Z}} \mathbb{Q}$. 
Is should be clear that $\varphi = \psi \circ \varphi _1$, as already pointed out in the introduction. 
We leave it as an (easy) exercise for the reader to check that these maps are well defined. 
This requires one to show that metabolic elements from any one group map to metabolic elements 
in the other groups. 
\vskip1mm

We conclude this section by reminding the reader of the isomorphism types of $\mathcal C_1$, $W(\mathbb{Q})$ and $\mathcal I(\mathbb Q)$ stated in the introduction: 

\centerline{
\xymatrix{
\mathcal C_1  \cong \mathbb{Z}^\infty \oplus \mathbb{Z}_2^\infty \oplus \mathbb{Z}_2^\infty  \ar[rr]^{\varphi_1} \ar[dr]_{\varphi }  &   & \mathcal I(\mathbb{Q}) \cong \mathbb{Z}^\infty \oplus \mathbb{Z}_2^\infty \oplus \mathbb{Z}_2^\infty   \ar[dl]^{\psi } \\
 &  W(\mathbb{Q})\mathcal  \cong \mathbb{Z} \oplus \mathbb{Z}_2^\infty \oplus \mathbb{Z}_2^\infty   &   \\
}
}
\vskip1mm
\noindent 
As already mentioned, Levine showed $\varphi _1$ to be injective. Clearly injectivity cannot hold for $\varphi$. 
However, given the above diagram, one cannot help but ask: \lq\lq How much loss of information is there
if one restricts $\varphi$ to the torsion subgroup of $\mathcal C_1$?\rq\rq \, 
As theorem \ref{main5} shows, the restriction of $\varphi$ to the torsion subgroup of $\mathcal C_1$ is unfortunately not injective. Nevertheless, examples \ref{example1} -- \ref{example4} show that 
$\varphi |_{Tor(\mathcal C_1)}$ contains significantly more information than just the knot determinant.  
\section{The linking matrices} \label{linking-section}
In this section we compute the linking matrix for $K=P(p_1,...,p_n)$  associated to a choice of oriented Seifert surface $\Sigma$ for $K$ along with a concrete basis for $H_1(\Sigma ;\mathbb{Z})$. The details of these computations for the three cases $(i)$--$(iii)$ from \eqref{categories} proceed in slightly different manners. 
\subsection{The case of $n,p_1,...,p_{n-1}$ odd and $p_n$ even} \label{one-one}
For the remainder of this subsection, we shall assume the conditions from its title with the additional constraints 
that $n\ge 3$ and $p_n\ne 0$.  
 
\begin{figure}[htb!] 
\centering
\includegraphics[width=14cm]{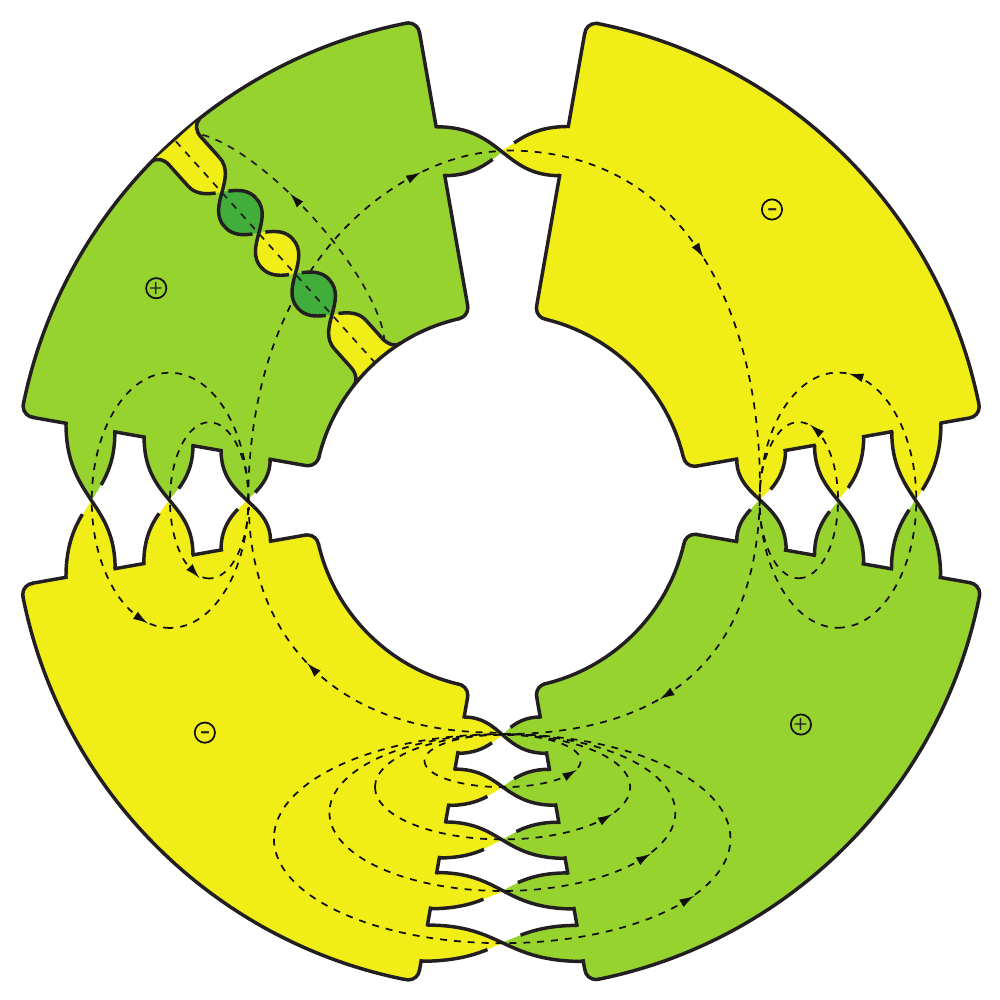}
\put(-50,248){$\alpha_2^2$}
\put(-70,228){$\alpha_1^2$}
\put(-173,83){$\alpha_1^3$}
\put(-158,65){$\alpha_2^3$}
\put(-143,48){$\alpha_3^3$}
\put(-128,30){$\alpha_4^3$}
\put(-335,228){$\alpha_1^4$}
\put(-364,240){$\alpha_2^4$}
\put(-300,340){$\delta$}
\put(-245,330){$\gamma$}
\put(-145,330){$\gamma$}
\put(-115,130){$\gamma$}
\put(-295,140){$\gamma$}
\caption{Our choice of Seifert surface $\Sigma_1$ for $P(p_1,...,p_n)$ for the case when $n, p_1, ..., p_{n-1}$ are odd and $p_n$ is even. Our example shows the knot $P(-1, 3, -5, 3, 4)$. The choices 
of generators for $H_1(\Sigma_1 ;\mathbb{Z})$ along with their orientations are indicated.}  \label{pic4}
\end{figure}
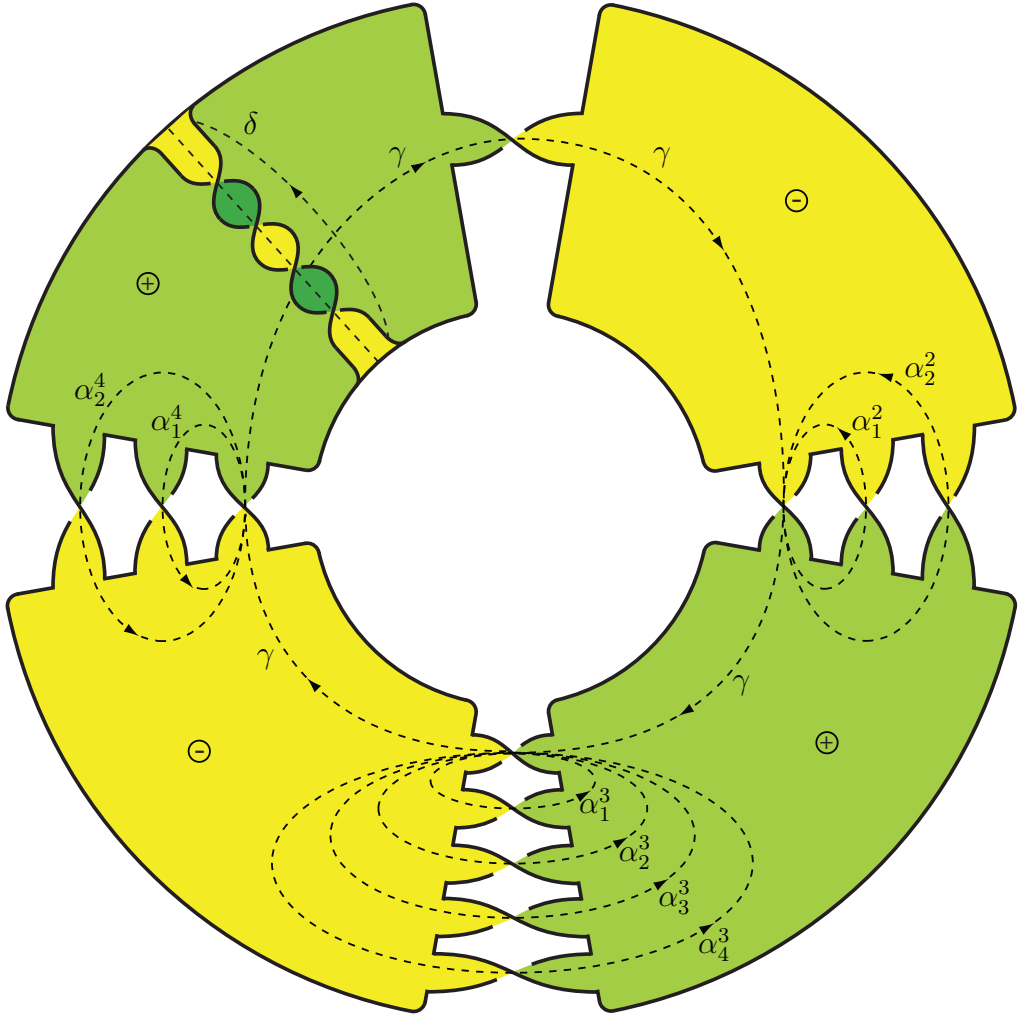

We start by recalling figure \ref{pic1} in which we chose a particular projection for the pretzel knot $P(p_1,...,p_n)$. 
We choose $\Sigma_1$ to be the Seifert surface for $K$ obtained from that projection via Seifert's algorithm (see for example \cite{rolfsen}). 
Specifically, $\Sigma_1$ consists of $n-1$ disks $D_1, ..., D_{n-1}$ of which $D_i$ and $D_{i+1}$ are connected 
with $|p_i|$ bands, each carrying a single half-twist whose handedness is determined by the sign of $p_i$ (in 
that the band obtains a right-handed twist if $p_i<0$ and a left-handed twist if $p_i>0$).  The disks $D_{n-1}$ and 
$D_1$ are similarly connected with $|p_{n-1}|$ bands. Finally, there is a band with $|p_n|$ 
half-twists (right-handed if $p_n>0$ and left-handed if $p_n<0$) both of whose ends are attached to $D_{1}$. 
Note that the genus of $\Sigma_1$ is $|p_1|+|p_2|+...+|p_{n-1}| + 3 -n $.
We label the bands connecting $D_i$ to $D_{i+1}$ 
by $B^i_1,...,B^i_{|p_i|}$ and we label those connecting $D_{n-1}$ to $D_1$ by $B^{n-1}_1,...,B^{n-1}_{|p_{n-1}|}$. 
The unique band with $|p_n|$ twists is labeled $B^n$. 
All of our conventions and labels are illustrated in figure \ref{pic4}.

With these preliminaries in place, we choose our basis 
\begin{equation} \label{basis1}
\mathcal B_1 =  \{ \alpha ^1_1,...,\alpha ^1_{|p_1| -1},\alpha ^2_1,...,\alpha ^2_{|p_2| -1},..., \alpha ^{n-1}_1,...,\alpha ^{n-1}_{|p_{n-1}| -1},\gamma, \delta \}
\end{equation}
for $H_1(\Sigma_1 ; \mathbb{Z})$ in the following way: 
\begin{enumerate}
\item We let $\alpha ^i_j$ to be the simple closed curve passing through the bands $B^i_1$ and$B^i_{j+1}$. 
\item We pick $\gamma$ to be the simple closed curve passing over the bands $B^1_1$, $B^2_1$, ..., $B^{n-1}_1$. 
\item The remaining curve $\delta$ passes once through the band $B^n$. 
\end{enumerate}
These curves, along with our orientation conventions, are also depicted in figure \ref{pic4}. The orientation 
of $\Sigma_1$ is determined by the normal vector field which points outwards from the page (and towards the 
reader) on all disks $D_1, D_3, D_5, ...$ and into the page (and away from the reader) on the disks 
$D_2, D_4, D_6, ...$ . These conventions are indicated by the symbols $\oplus$ and $\ominus$ respectively in figure \ref{pic4}.

With these definitions in place, we are ready to start computing entries in the linking matrix $\mathcal L = [\mathcal \ell _{i,j}]$
where $\ell _{i,j} = \ell k (x_i , x_j)$. Here $x_i$ is the $i$-th element of the basis $\mathcal B_1$ and 
$\ell k (x_i , x_j)$ is the linking number of $x_i$ and $x_j^+$. The latter is a small push-off of $x_j$ in the direction of the normal vector field on $\Sigma_1$ determined by its orientation.  

Seeing as the loops $\alpha^i _k$ and $\alpha^j_m$ are disjoint for any choice of $i\ne j$, we find that 
$\ell k(\alpha ^i_k, \alpha^j_m) = \ell k(\alpha^j_m,\alpha ^i_k) =0$ for any choices of $i,j,k,m$ with $i\ne j$. 
For the same reason, one also obtains $\ell k (\alpha ^i_k, \delta) = \ell k (\delta, \alpha ^i_k) = 0$ for any choices
of $i,k$. 

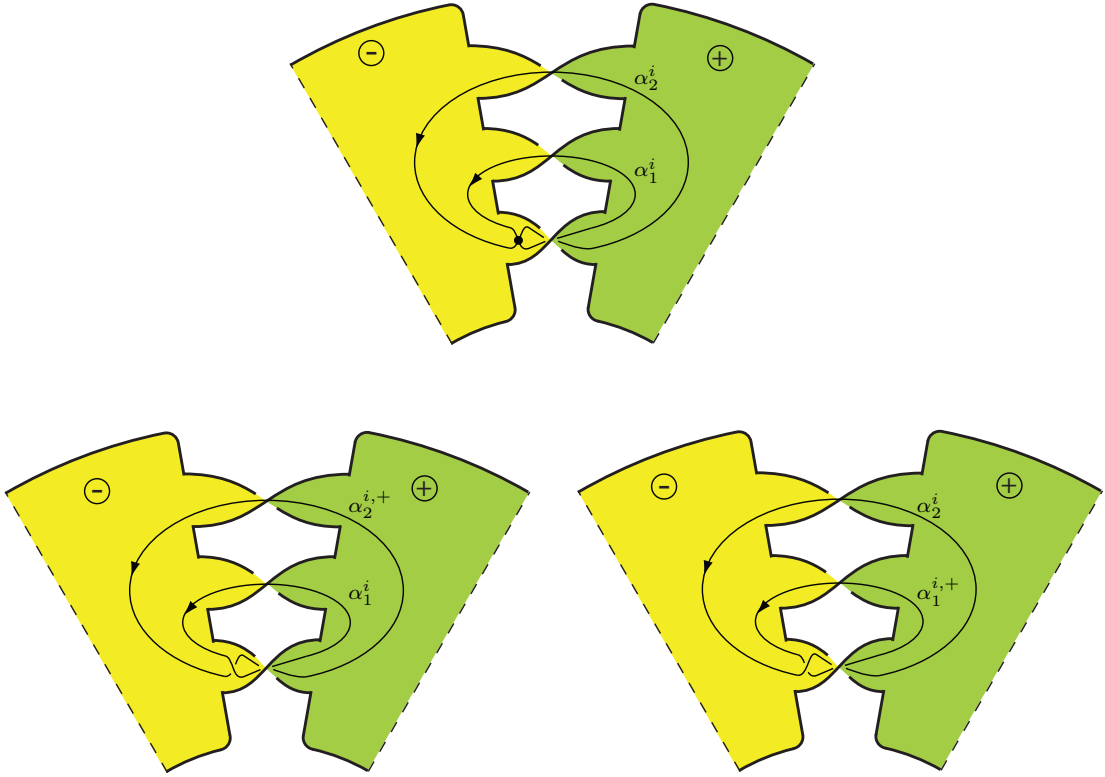
\begin{figure}[htb!] 
\centering
\includegraphics[width=15cm]{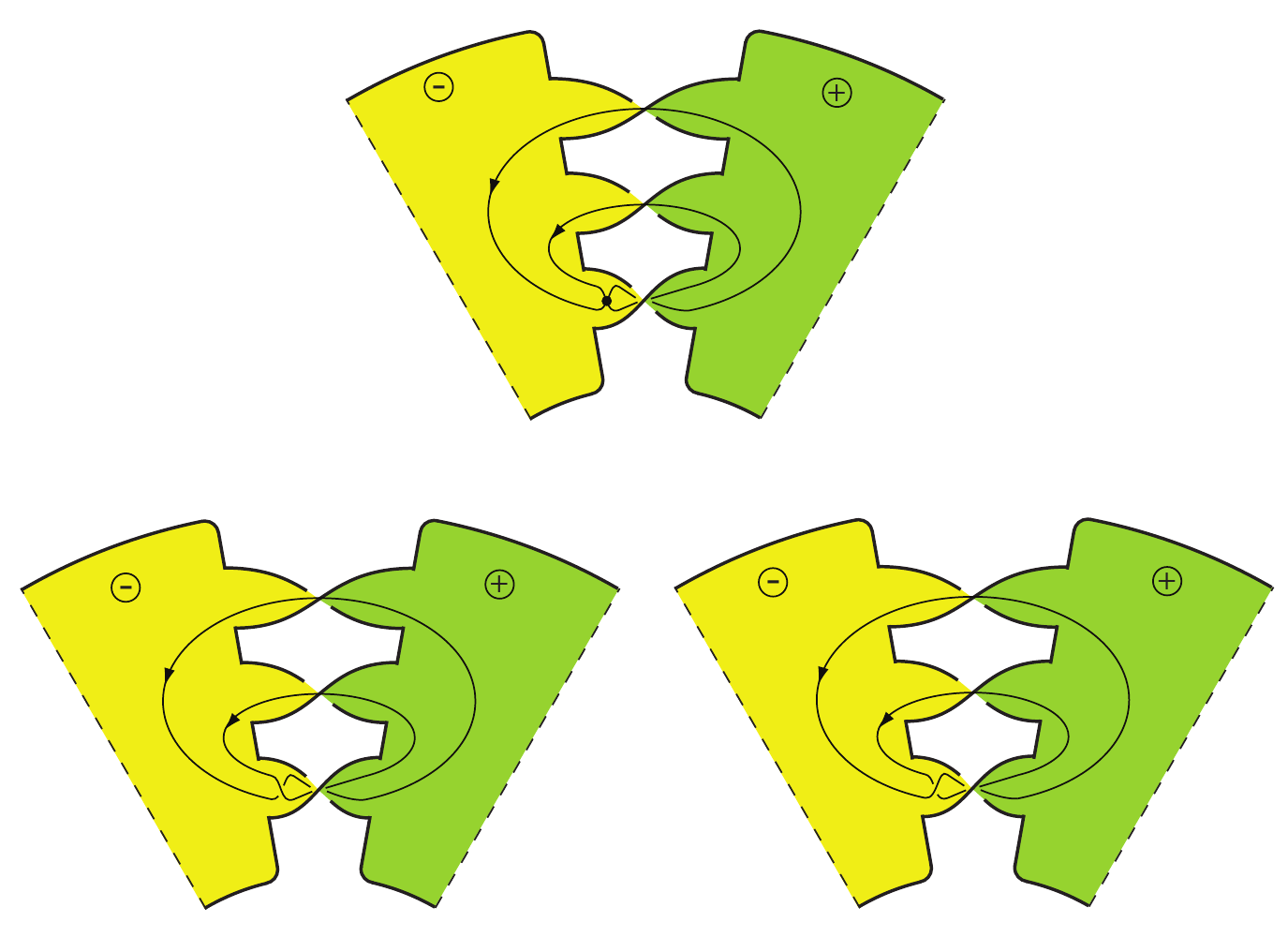}
\put(-182,235){\tiny $ \alpha_1^i$}
\put(-182,270){\tiny $\alpha_2^i$}
\put(-290,75){\tiny $ \alpha_1^i$}
\put(-290,108){\tiny $\alpha_2^{i,+}$}
\put(-75,75){\tiny $ \alpha_1^{i,+}$}
\put(-75,108){\tiny $\alpha_2^i$}
\caption{This figure computes the linking $\ell k (\alpha ^i_1, \alpha ^i_2)$ when $i$ is even. The two push-offs 
$\alpha ^{i,+}_1$ and $\alpha ^{i,+}_2$ of $\alpha ^i_1$ and $\alpha ^i_2$ respectively, are shown in the bottom two pictures. The linking of the two is readily computed from these.}  \label{pic5}
\end{figure}

The contribution of the subset $\{\alpha ^i_1,...,\alpha ^i_{|p_i|-1}\}$ of $\mathcal B$ to the linking form $\mathcal L$, only depends on $p_i$. To see how, let us introduce the $n\times n$ matrices $X_n$  and $Y_n = X_n+X_n^\tau$ by the formulae
\begin{equation} \label{X-n}
 X_n = \left[ 
\begin{array}{rrrrrr}
1 & 0 &0 & ... & 0 & 0 \cr
1 & 1 & 0 &  ... & 0 & 0 \cr
1 & 1 & 1 &  ... & 0 & 0 \cr
\vdots & \vdots & \vdots   & \ddots & \vdots & \vdots  \cr 
1 & 1 & 1 &  ... & 1 & 0  \cr
1 & 1 & 1 &  ... & 1 & 1 \cr
\end{array}
\right] \quad \quad \mbox{ and } \quad \quad 
Y_n =  \left[ 
\begin{array}{rrrrrr}
2 & 1 &1 & ... & 1 & 1 \cr
1 & 2 & 1 &  ... & 1 & 1 \cr
1 & 1 & 2 &  ... & 1 & 1 \cr
\vdots & \vdots  & \vdots  & \ddots & \vdots & \vdots  \cr 
1 & 1 & 1 &  ... & 2 & 1  \cr
1 & 1 & 1 &  ... & 1 & 2 \cr
\end{array}
\right]
\end{equation}
By consulting figure \ref{pic4}, one finds that  
\begin{equation} \label{linkingeq-1}
 \begin{array}{ll}
\ell k(\alpha ^i_k , \alpha ^i_m) = \left\{ 
\begin{array}{rcl}
0 & \quad ; \quad & k<m \cr
-1 & \quad ; \quad &  k\ge m 
\end{array}
\right. & \quad \quad \quad \mbox{ if } p_i >0 \mbox{ and } i\mbox{ is even.}\cr
 & \cr 
\ell k(\alpha ^i_k , \alpha ^i_m) = \left\{ 
\begin{array}{rcl}
\phantom{-}1 & \quad ; \quad & k\le m \cr
0 & \quad ; \quad &  k> m 
\end{array}
\right. & \quad \quad \quad \mbox{ if } p_i <0 \mbox{ and } i\mbox{ is even.}\cr
& \cr 
\ell k(\alpha ^i_k , \alpha ^i_m) = \left\{ 
\begin{array}{rcl}
 -1 & \quad ; \quad & k\le m \cr
 0 & \quad ; \quad &  k> m 
\end{array}
\right. & \quad \quad \quad \mbox{ if } p_i >0 \mbox{ and } i\mbox{ is odd.}\cr
 & \cr 
\ell k(\alpha ^i_k , \alpha ^i_m) = \left\{ 
\begin{array}{rcl}
 0 & \quad ; \quad & k <  m \cr
 \phantom{-}1& \quad ; \quad &  k\ge m 
\end{array}
\right. & \quad \quad \quad \mbox{ if } p_i <0 \mbox{ and } i\mbox{ is odd.}
\end{array}
\end{equation}
The case of $p_i>0$ and $i$ even is singled out in figure \ref{pic5}. 
From this we find that $\mathcal L_i$, the restriction of the linking form $\mathcal L$ to the 
$Span (  \alpha ^i_1, ..., \alpha ^i_{|p_i|-1})$, with respect to the basis $\{ \alpha ^i_1,...,\alpha ^i_{|p_i|-1}\}$ takes on one of 4 possible forms:
$$ \mathcal{L} _i = \left\{
\begin{array}{rl}
-X_{|p_i|-1} & \quad ; \quad \mbox{ if } p_i >0 \mbox{ and } i\mbox{ is even.}   \cr & \cr 
X^\tau _{|p_i|-1} & \quad ; \quad \mbox{ if } p_i <0 \mbox{ and } i\mbox{ is even.}   \cr & \cr
-X^\tau _{|p_i|-1} & \quad ; \quad \mbox{ if } p_i >0 \mbox{ and } i\mbox{ is odd.}   \cr & \cr 
X _{|p_i|-1} & \quad ; \quad \mbox{ if } p_i <0 \mbox{ and } i\mbox{ is odd.}  
\end{array}
\right. 
$$
\vskip2mm
In each of the four cases above, the matrix representing $\mathcal L_i + \mathcal L_i^\tau$ can then be expressed as 
\begin{equation} 
\mathcal L_i + \mathcal L_i^\tau  = - Sign(p_i) \, \, Y_{|p_i|-1} 
\end{equation}
\vskip2mm

Having worked out all of the linking numbers $\ell k (\alpha ^i_k, \alpha ^j _m)$, we now turn to exploring how 
$\gamma$ and $\delta$ contribute to $\mathcal L$. Their linking numbers with the various other 
curves from the basis $\mathcal B$ are easily read off from figure \ref{pic4}:
\begin{align} \label{linkingeq-2}
\ell k (\gamma, \gamma ) & = -(Sign(p_1) + Sign (p_2) + ... + Sign (p_{n-1}))/2 \cr
\ell k (\gamma, \delta) & = 0 \cr 
\ell k (\delta, \gamma ) & = 1 \cr
\ell k (\delta, \delta) & = p_n/2 
\end{align}
while the linking numbers of $\gamma$ with the various $\alpha ^i_k$ are 
\begin{equation} \label{linkingeq-3}
\begin{array}{rl}
\ell k (\gamma, \alpha ^i_k) & = \left\{
\begin{array}{rr} 
-1 &  \quad ; \quad \mbox{ if } p_i >0 \mbox{ and } i\mbox{ is even.}   \cr
0 & \quad ; \quad \mbox{ if } p_i <0 \mbox{ and } i\mbox{ is even.}   \cr 
0 & \quad ; \quad \mbox{ if } p_i >0 \mbox{ and } i\mbox{ is odd.}   \cr 
 1 & \quad ; \quad \mbox{ if } p_i <0 \mbox{ and } i\mbox{ is odd.}  
\end{array}
\right. \cr & \cr 
\ell k (\alpha ^i_k , \gamma) & = \left\{
\begin{array}{rr} 
0 &  \quad ; \quad \mbox{ if } p_i >0 \mbox{ and } i\mbox{ is even.}   \cr
1 & \quad ; \quad \mbox{ if } p_i <0 \mbox{ and } i\mbox{ is even.}   \cr 
-1 & \quad ; \quad \mbox{ if } p_i >0 \mbox{ and } i\mbox{ is odd.}   \cr 
 0 & \quad ; \quad \mbox{ if } p_i <0 \mbox{ and } i\mbox{ is odd.}  
\end{array}
\right. 
\end{array}
\end{equation}
As earlier, we see that while $\ell k (\alpha ^i_k , \gamma)$ and $\ell k (\gamma, \alpha ^i_k)$ depend on a number
of cases, the quantity $\ell k (\gamma, \alpha ^i_k) + \ell k (\alpha ^i_k , \gamma)$ always equals $ -Sign (p_i)$.
We are thus in a position to assemble all the pieces. 
\begin{theorem} \label{linkingform1}
Let $n,p_1,...,p_{n-1}$ be odd integers with $n \ge 3$ and let $p_n \ne 0$ be an even integer. 
To keep notation below at bay, let us also introduce the abbreviations 
$$ s_i = -Sign(p_i) \quad \quad  \quad s=s_1+...+s_{n-1}  \quad \quad  \quad   \rho_i=|p_i|-1$$
Then the symmetrized linking form $\mathcal L + \mathcal L^\tau $ of the pretzel knot $P(p_1,...,p_n)$ associated to the oriented Seifert surface $\Sigma_1$ and the basis   
$$ \mathcal B =  \{ \alpha ^1_1,...,\alpha ^1_{|p_1| -1},\alpha ^2_1,...,\alpha ^2_{|p_2| -1},..., \alpha ^{n-1}_1,...,\alpha ^{n-1}_{|p_{n-1}| -1},\gamma, \delta \}$$
of $H_1(\Sigma_1 ; \mathbb{Z})$ as chosen above (see specifically figure \ref{pic4}), has the form
$$\mathcal L + \mathcal L^\tau  = \left[
\begin{array}{ccc|ccc|c|ccc|c|c}
   &   &   &   &   &    &  \dots   &    &       &     &  s_1 & 0  \cr 
   & s_1 Y_{\rho_1} & &   & 0  &     &  \dots      &    &  0  &     &  \vdots & \vdots  \cr 
   &  & &   &   &    &   \dots  &    &   &         &  s_1 & 0  \cr \hline 
   &  & &   &   &    &   \dots  &    &   &         &  s_2 & 0  \cr
   & 0 & &   & s_2Y_{\rho_2}  &    &  \dots      &    &  0  &     &  \vdots & \vdots  \cr
   &  & &   &   &    &   \dots  &    &   &         &  s_2 & 0  \cr \hline 
  \vdots  &  \vdots &  \vdots &  \vdots  &  \vdots  &  \vdots   &  \ddots  &    &   &         &  \vdots  &   \vdots \cr \hline
   &  & &   &   &    &   &    &   &         & s_{n-1}  & 0   \cr 
   & 0 & &   & 0  &    &   &    & s_{n-1}Y_{\rho_{n-1}}  &         & \vdots  & \vdots   \cr 
   &  & &   &   &    &   &    &   &         & s_{n-1}  & 0   \cr  \hline 
  s_1& \dots &s_1 & s_2  & \dots   & s_2   & \dots  &   s_{n-1} & \dots  & s_{n-1}        & s & 1   \cr \hline 
   0 & \dots  & 0 & 0  & \dots   & 0    & \dots  &   0 &  \dots  &    0     & 1 & p_n   \cr 
\end{array} 
\right]
$$ 
The matrices $Y_\rho$ are as introduced in \eqref{X-n}. 
\end{theorem}
%
%
%
\subsection{The case of $n$ even, $p_1,...,p_{n-1}$ odd and $p_n$ even}
We turn to the next case of choice of parities of $n, p_1,...,p_n$ and pick it for the remainder of this 
section to be as listed in the title. We also keep our additional assumptions of $n\ge 3$ and $p_n\ne 0$. 

\begin{figure}[htb!] 
\centering
\includegraphics[width=15cm]{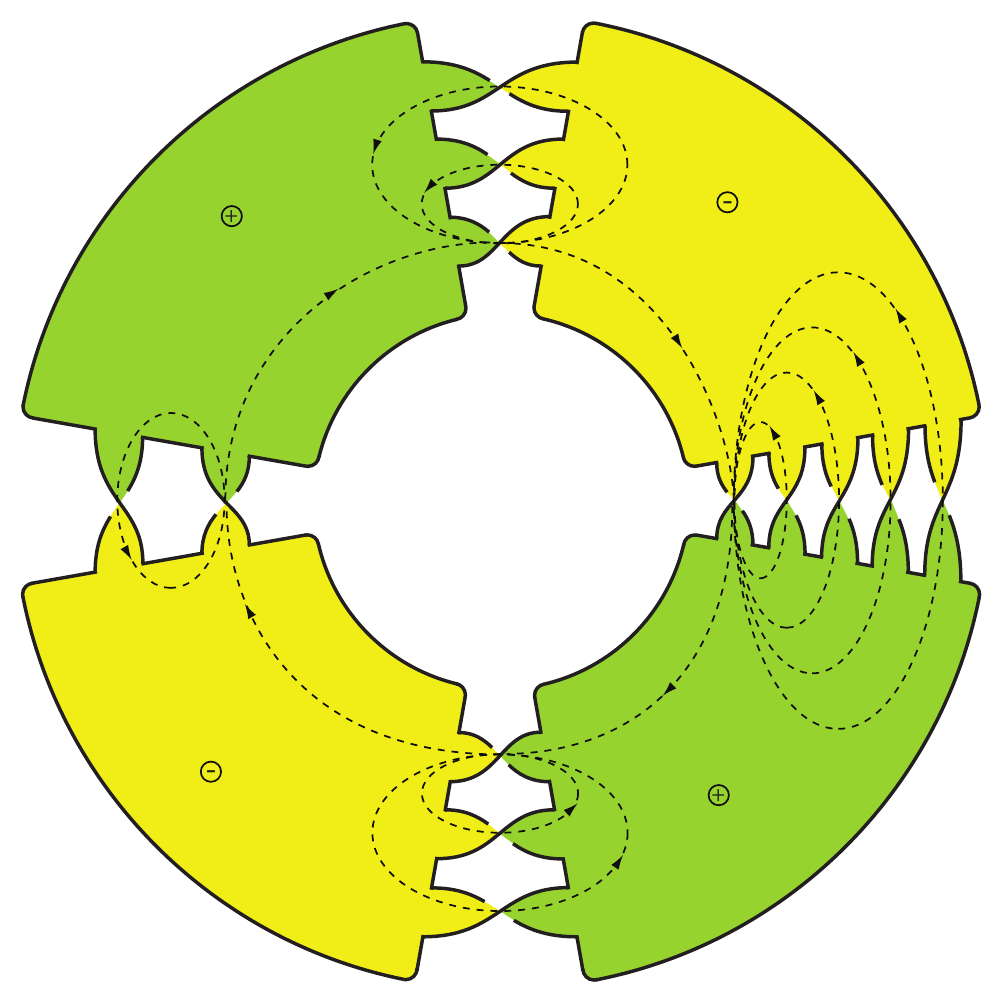}
\put(-183,350){$\alpha_1^1$}
\put(-175,383){$\alpha_2^1$}
\put(-100,247){$\alpha_1^2$}
\put(-82,265){$\alpha_2^2$}
\put(-65,282){$\alpha_3^2$}
\put(-47,301){$\alpha_4^2$}
\put(-183,72){$\alpha_1^3$}
\put(-175,40){$\alpha_2^3$}
\put(-355,255){$\alpha_1^4$}
\put(-300,280){$\gamma$}
\put(-160,288){$\gamma$}
\put(-160,130){$\gamma$}
\put(-285,130){$\gamma$}
\caption{Our choice of Seifert surface $\Sigma_2$ for $P(p_1,...,p_n)$ for the case when $n$ is even, $p_1, ..., p_{n-1}$ are odd and $p_n$ is even. Our example shows the knot $P(3,-5, 3, 2)$. The choices 
of generators for $H_1(\Sigma_2 ;\mathbb{Z})$ with their orientations are indicated.}  \label{pic7}
\end{figure}

The Seifert surface $\Sigma_2$ that we choose for $P(p_1,...,p_n)$
and the  preferred basis $\mathcal B_2$ for $H_1(\Sigma_2 ; \mathbb{Z})$ are very much like in the case considered in section \ref{one-one}.  Specifically, we let $\Sigma_2$ be obtained from $\Sigma_1$ ($\Sigma_1$ is the Seifert surface from section \ref{one-one}) by simply deleting its unique band with and even number of half-twists and allowing the number of bands which connect the disks $D_{n}$ and $D_1$ to be an even number, namley $|p_n|$. We then arrive 
at a surface $\Sigma_2$ as in figure \ref{pic7}. The same figure also indicates  our choice of basis 
$$ \mathcal B_2 =  \{ \alpha ^1_1,...,\alpha ^1_{|p_1| -1},\alpha ^2_1,...,\alpha ^2_{|p_2| -1},..., \alpha ^{n}_1,...,\alpha ^{n}_{|p_{n}| -1},\gamma \}$$
for $H_1(\Sigma _2;\mathbb{Z})$ which is identical to $\mathcal B_1$ from \eqref{basis1} safe that we are presently no longer requiring the generator $\delta$. The orientation convention is as in the previous section and is again indicated by a $\oplus$ 
and $\ominus$ in figure \ref{pic7}.  

The linking numbers between the various $\alpha ^i_k$ and $\alpha ^j_m$ and indeed between the $\alpha ^i_k$ 
and $\gamma$ are identical to those found in section \ref{one-one}. We thus immediately arrive at the 
analogue of theorem \ref{linkingform1}:
\begin{theorem} \label{linkingform2}
Let $n\ge 3$ be an even integer and let $p_1,...,p_{n-1}$ be odd integers and $p_n\ne 0$ an even integer. 
Let us re-introduce the abbreviations 
$$ s_i = -Sign(p_i) \quad \quad  \quad s=s_1+...+s_{n-1}  \quad \quad  \quad   \rho_i=|p_i|-1$$
Then the symmetrized linking form $\mathcal L + \mathcal L^\tau $ of the pretzel knot $P(p_1,...,p_n)$ associated to the oriented Seifert surface $\Sigma_2$ and the basis   
$$ \mathcal B_2 =  \{ \alpha ^1_1,...,\alpha ^1_{|p_1| -1},\alpha ^2_1,...,\alpha ^2_{|p_2| -1},..., \alpha ^{n}_1,...,\alpha ^{n}_{|p_{n}| -1},\gamma \}$$
of $H_1(\Sigma_2 ; \mathbb{Z})$ as chosen above (see specifically figure \ref{pic7}) takes the form
$$\mathcal L + \mathcal L^\tau  = \left[
\begin{array}{ccc|ccc|c|ccc|c}
   &   &   &   &   &    &  \dots   &    &       &     &  s_1   \cr 
   & s_1 Y_{\rho_1} & &   & 0  &     &  \dots      &    &  0  &     &  \vdots   \cr 
   &  & &   &   &    &   \dots  &    &   &         &  s_1   \cr \hline 
   &  & &   &   &    &   \dots  &    &   &         &  s_2   \cr
   & 0 & &   & s_2Y_{\rho_2}  &    &  \dots      &    &  0  &     &  \vdots   \cr
   &  & &   &   &    &   \dots  &    &   &         &  s_2   \cr \hline 
  \vdots  &  \vdots &  \vdots &  \vdots  &  \vdots  &  \vdots   &  \ddots  &    &   &         &  \vdots   \cr \hline
   &  & &   &   &    &   &    &   &         & s_{n}     \cr 
   & 0 & &   & 0  &    &   &    & s_{n}Y_{\rho_{n}}  &         & \vdots     \cr 
   &  & &   &   &    &   &    &   &         & s_{n}     \cr  \hline 
  s_1& \dots &s_1 & s_2  & \dots   & s_2   & \dots  &   s_{n} & \dots  & s_{n}        & s 
\end{array} 
\right]
$$ 
The matrices $Y_\rho$ are as in  \eqref{X-n}. 
\end{theorem}
%
%
%
\subsection{The case of $n$ and $p_1,...,p_{n}$ odd}
\begin{figure}[htb!] 
\centering
\includegraphics[width=15cm]{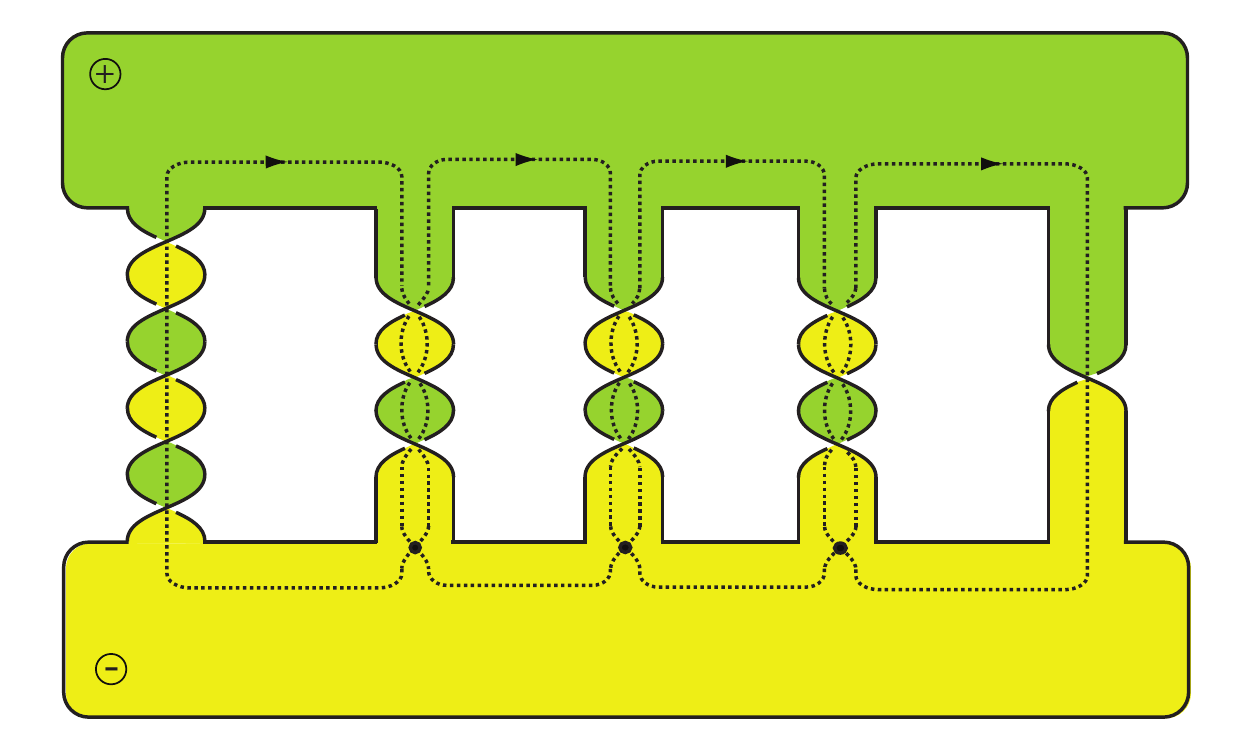}
\put(-332,210){$\alpha_1$}
\put(-255,210){$\alpha_2$}
\put(-180,210){$\alpha_3$}
\put(-100,210){$\alpha_4$}
\caption{Our choice of Seifert surface $\Sigma _3$ for $P(p_1,...,p_n)$ for the case when $n$ and $p_1, ..., p_{n}$ are odd. Our example shows the knot $P(5,-3, 3, -3,-1)$. The choices 
of generators for $H_1(\Sigma_3 ;\mathbb{Z})$ with their orientations are indicated.}  \label{pic6}
\end{figure}

In this section we consider the remaining case where all of $n, p_1,...,p_n$ are odd with $n\ge 3$. We start 
by picking a Seifert surface $\Sigma _3$ for $P(p_1,...,p_n)$ which is this time obtained by taking two 
disks and connecting them by $n$ bands $B_1$, ... $B_n$, each with $|p_i|$ half twists (right-handed twists if $p_i>0$ and 
left-handed twists if $p_i<0$). The thus obtained surface looks as in figure \ref{pic6}.
We next choose a basis 
$$\mathcal B_3 =\{ \alpha _1, ...,\alpha _{n-1} \}$$
of $H_1(\Sigma _3 ; \mathbb{Z})$ by letting $\alpha_i$ be the curve on $\Sigma _3$ which runs through the 
bands $B_i$ and $B_{i+1}$. The orientation conventions for the $\alpha _i$ and indeed the orientation 
for $\Sigma _2$ itself (indicated again by a $\oplus$ and a $\ominus$) are depicted in figure \ref{pic6}. 

The linking form in this basis is rather easy to determine. Note first that $\ell k (\alpha _i , \alpha _j) = 0$ whenever
$|i-j|\ge 2$. On the other hand, by inspection from figure \ref{pic6}, it follows that 
\begin{align} \nonumber
\ell k (\alpha _i, \alpha _{i} )  = \frac{p_i+p_{i+1}}{2} \quad \quad \ell k (\alpha _i , \alpha _{i+1}) = -\frac{p_{i+1}+1}{2}
\quad \quad \ell k (\alpha _{i+1} , \alpha _{i}) = -\frac{p_{i+1}-1}{2}
\end{align}
With this in place, here is the analogue of theorems \ref{linkingform1} and \ref{linkingform2} for the present case. 
\begin{theorem} \label{linkingform3}
Let $n\ge 3$ be an odd integer and let $p_1,...,p_{n}$ be any odd integers.  
Then the symmetrized linking form $\mathcal L + \mathcal L^\tau $ of the pretzel knot $P(p_1,...,p_n)$ associated to the oriented Seifert surface $\Sigma_3$ and the basis $ \mathcal B_3 =  \{ \alpha _1,...,\alpha _{n -1} \}$
of $H_1(\Sigma_3 ; \mathbb{Z})$ as chosen above (see figure \ref{pic6}) takes the form
$$\mathcal L + \mathcal L^\tau  = \left[
\begin{array}{cccccccc}
p_1+p_2 & -p_2 & 0 &0&  \dots  & 0 & 0 & 0\cr  
-p_2 & p_2 + p_3& -p_3 & 0&  \dots  & 0 & 0 & 0 \cr  
0& -p_3 & p_3 + p_4& -p_4 & \dots &  0  & 0  & 0 \cr  
\vdots  &  &  & \vdots &    &  & &  \vdots  \cr  
0& 0 & 0 & 0 & \dots &  -p_{n-2}  & p_{n-2}+p_{n-1}  & -p_{n-1} \cr  
0& 0 & 0 & 0 & \dots &  0  & -p_{n-1}  & p_{n-1}+p_{n} 
\end{array} 
\right]
$$ 
\end{theorem}
%
%
%
\section{Diagonalizing the linking matrices} \label{diagonalizing-section}
In this section we show how one can diagonalize the matrices $\mathcal L + \mathcal L^\tau$ obtained in theorems \ref{linkingform1}, \ref{linkingform2} and \ref{linkingform3}. We do this essentially using the Gram-Schmidt process on $(\langle \cdot , \cdot \rangle,H_1(\Sigma ; \mathbb{Q}))$ with $\langle x, y \rangle = \ell k (x,y) + \ell k (y,x)$.  We need to exercise a bit of care since, while $\langle \cdot , \cdot \rangle$ is non-degenerate, it is by no means definite and square zero vectors do exist. 

Once $\mathcal L + \mathcal L^\tau$ has been diagonalized, it is an easy matter to read off the rational Witt class of $\mathcal L + \mathcal L^\tau$ in terms of the generators of $W(\mathbb{Q})$. 
\subsection{The Gram-Schmidt procedure and reduction} \label{GSRed}
We start by reminding the reader of the Gram-Schmidt process on an arbitrary finite dimensional inner product 
space $(\langle \cdot , \cdot \rangle,V)$. By convention, such an inner product $\langle \cdot , \cdot \rangle$ is 
assumed to be positive definite. We then address the issue of square zero vectors in 
$(\langle \cdot , \cdot \rangle,H_1(\Sigma ; \mathbb{Q}))$. 
\begin{theorem}[Gram-Schmidt] \label{gramschmidttheorem}
Let $\{f_1,...,f_n\}$ be a basis for the inner product space $(\langle \cdot , \cdot \rangle,V)$ and let $\{e_1,...,e_n\}$
be the set of vectors obtained as 
\begin{align} \nonumber
e_1 & = f_1 \cr
e_2 & = f_2 - \frac{\langle f_2, e_1\rangle}{\langle e_1, e_1\rangle} e_1 \cr 
e_3 & = f_3 - \frac{\langle f_3, e_2\rangle}{\langle e_2, e_2\rangle} e_2 - \frac{\langle f_3, e_1\rangle}{\langle e_1, e_1\rangle} e_1 \cr 
& \vdots \cr 
e_n & = f_n - \frac{\langle f_n, e_{n-1}\rangle}{\langle e_{n-1}, e_{n-1}\rangle} e_{n-1} - ... - \frac{\langle f_n, e_1\rangle}{\langle e_1, e_1\rangle} e_1 \cr 
\end{align}
Then $\{e_1,...,e_n\}$ is an orthogonal basis for $V$ and $Span\{e_1,...,e_i\} = Span\{f_1,...,f_i\}$ 
for each $i\le n$. 
\end{theorem}
\begin{remark} \label{gramgcd}
In order to keep the scalars in our computations integral (rather than rational and non-integral), we will often use the 
slightly modified Gram-Schmidt process by which we set 
$$ e_i  = d_i \cdot \left( f_i - \frac{\langle f_i, e_{i-1}\rangle}{\langle e_{i-1}, e_{i-1}\rangle} e_{i-1} - ... - \frac{\langle f_i, e_1\rangle}{\langle e_1, e_1\rangle} e_1\right)$$
where $d_i$ is some common multiple of $\langle e_{1}, e_{1}\rangle, ... ,\langle e_{i-1}, e_{i-1}\rangle$. 
Clearly, the thus created set $\{e_1,...,e_n\}$ is still an orthogonal basis for any choice of $d_i\ne 0$. 
\end{remark}
The next theorem addresses the failure of the Gram-Schmidt procedure in the presence of square zero vectors
(on non-definite inner product spaces). 
The result should be viewed as an iterative prescription to be applied as many times in the Gram-Schmidt process as is the number of square zero vectors $e_i$ encountered. 
\begin{theorem} \label{squarezero}
Let $(\lr, V)$ be a pair consisting of a finite dimensional $\mathbb{F}$-vector space $V$ and a non-degenerate bilinear symmetric form $\lr$. Let $\{f_1,...,f_n\}$ be a basis for $V$ and let, for some $m<n$, $\{e_1,...,e_m\}$ 
be obtained from $\{f_1,...,f_n\}$ as in theorem \ref{gramschmidttheorem} (or alternatively as in remark \ref{gramgcd}). Assume that $\langle e_i ,e_i \rangle \ne 0$ for all $i< m$ but that $\langle e_m ,e_m \rangle = 0$. Additionally, suppose also that $\langle e_m,f_{m+1}\rangle \ne 0$ (which can always be achieved by a simple 
reordering, if necessary, of $f_{m+1},...,f_n$). 

Then $(\lr, V)$ is equal to $(\lr',V')$ in the Witt ring $W(\mathbb{F})$ where 
$$V' =Span(e_1,...,e_{m-1},f'_{m+2},...,f'_n) \quad  \quad \mbox{ and } \quad \quad \lr' = \lr |_{V'\times V'} 
$$
with
\begin{align} \nonumber
f'_{m+1} & = f_{m+1} - \sum _{j=1}^{m-1} \frac{\langle f_{m+1},e_j \rangle }{\langle e_{j},e_{j}\rangle} e_j   \cr
f''_{m+k} & = f_{m+k} - \sum _{j=1}^{m-1} \frac{\langle f_{m+k},e_j\rangle }{\langle e_j , e_j \rangle } e_j \cr
f'_{m+k} & = f''_{m+k} -  \frac{\langle f''_{m+k},e_{m}\rangle }{\langle f'_{m+1},e_m\rangle}f'_{m+1} 
-\frac{\langle f''_{m+k} , f'_{m+1}\rangle \cdot \langle f'_{m+1} , e_m \rangle  - \langle f''_{m+k},e_m\rangle\cdot \langle f'_{m+1},f'_{m+1}\rangle }{\langle f'_{m+1},e_m\rangle \cdot \langle f'_{m+1},e_m\rangle } e_m
\end{align}
where the last two equations are valid for $k\ge 2$. 
\end{theorem}
\begin{proof}
Let $A$ be the symmetric non-degenerate $n\times n$ matrix representing $\lr$ with respect to the basis $\{e_1,...,e_{m-1}\}\cup \{e_m,f_{m+1},...,f_n\}$.  Then $A$ is of the form 
$$A=  \left[
\begin{array}{ccc|cccc}
\langle e_1,e_1 \rangle  &\dots  &0 & 0  & \langle e_1,f_{m+1} \rangle  &...& \langle e_1,f_{n} \rangle \cr
 \vdots & \ddots & \vdots   & \vdots & \vdots & \ddots & \vdots  \cr
 0 & \dots &\langle e_{m-1},e_{m-1} \rangle   & 0 & \langle e_{m-1},f_{m+1} \rangle & \dots & \langle e_{m-1},f_{m+1} \rangle \cr \hline 
0  & \dots &0 &0 &   \langle e_{m},f_{m+1} \rangle &...& \langle e_{m},f_{n} \rangle \cr
 \langle f_{m+1},e_{1} \rangle & \dots  &  \langle f_{m+1},e_{m-1} \rangle &   \langle f_{m+1},e_m \rangle &   \langle f_{m+1},f_{m+1} \rangle&  ...& \langle f_{m+1},f_{n} \rangle \cr
 \vdots & \ddots  & \vdots & \vdots & \vdots &\ddots  & \vdots \cr
  \langle f_{n},e_{1} \rangle & \dots  &  \langle f_{n},e_{m-1} \rangle &   \langle f_{n},e_m \rangle &   \langle f_{n},f_{m+1} \rangle&  ...& \langle f_{n},f_{n} \rangle
\end{array}
\right]
$$
For $k\ge 1$, let $f''_{m+k}$ be given by 
$$ f''_{m+k} = f_{m+k} - \sum _{j=1}^{m-1} \frac{\langle f_{m+k},e_j\rangle }{\langle e_j , e_j \rangle } e_j $$
so that $\langle f''_{m+k}, e_i\rangle  = 0$ for all $k\ge 1$ and all $i\le m-1$. Thus the matrix $A''$ representing
$\mathcal L + \mathcal L^\tau$ with respect to the basis $\{ e_1,...e_m,f''_{m+1},...,f''_{n}\}$ looks like 
$$A '' =  \left[
\begin{array}{ccc}
\langle e_1,e_1 \rangle  &\dots  &0  \cr
\vdots & \ddots & \vdots  \cr
0 & \dots &\langle e_{m-1},e_{m-1} \rangle 
\end{array} \right] \oplus
\left[
\begin{array}{cccc} 
0 &   \langle e_{m},f_{m+1} \rangle &...& \langle e_{m},f_{n} \rangle \cr
\langle f_{m+1},e_m \rangle &   \langle f''_{m+1},f''_{m+1} \rangle&  ...& \langle f''_{m+1},f''_{n} \rangle \cr
\vdots & \vdots &\ddots  & \vdots \cr
\langle f_{n},e_m \rangle &   \langle f''_{n},f''_{m+1} \rangle&  ...& \langle f''_{n},f''_{n} \rangle
\end{array}
\right]
$$
Note that $\langle e_m,f''_{m+k} \rangle = \langle e_m,f_{m+k} \rangle $ for all $k\ge 1$. 
To simplify the second summand, we introduce a further change of basis by setting 
$$f'_{m+k} = f''_{m+k} -  \frac{\langle f''_{m+k},e_{m}\rangle }{\langle f''_{m+1},e_m\rangle}f''_{m+1} 
-\frac{\langle f''_{m+k} , f''_{m+1}\rangle \cdot \langle f''_{m+1} , e_m \rangle  - \langle f''_{m+k},e_m\rangle\cdot \langle f''_{m+1},f''_{m+1}\rangle }{\langle f''_{m+1},e_m\rangle \cdot \langle f''_{m+1},e_m\rangle } e_m
$$
for all $k\ge 2$ and for convenience, set $f'_{m+1}=f''_{m+1}$. 
A quick check reveals that now 
$$ \langle f'_{m+k} , e_m \rangle =0 \quad \quad \mbox{  and } \quad \quad \langle f'_{m+k} , f'_{m+1} \rangle =0 \quad 
\quad \forall k\ge 2$$
Therefore the second summand of $A''$ above, when expressed with respect to the basis $\{e_m,f'_{m+1},...,f'_n\}$, takes the form 
$$ \left[
\begin{array}{cc}
0 & \langle e_m,f_{m+1} \rangle \cr
\langle f_{m+1},e_m \rangle &   \langle f'_{m+1},f'_{m+1}  \rangle 
\end{array}
\right] \oplus 
\left[
\begin{array}{ccc}
 \langle f'_{m+2},f'_{m+2} \rangle  & ... & \langle f'_{m+2},f'_{n} \rangle \cr
\vdots &  & \vdots \cr
 \langle f'_{n},f'_{m+2} \rangle &...& \langle f'_{n},f'_{n} \rangle 
\end{array}
\right]
$$
Since the first summand is metabolic and therefore equals zero in $W(\mathbb{F})$, the claim of the theorem 
follows. 
\end{proof}
We shall refer to the passage from $(\lr,V)$ to $(\lr',V')$, as described in theorem \ref{squarezero}, as {\em reduction}, seeing as the dimension of $V$ gets reduced by $2$ in the process.  
\subsection{The case of $n, p_1,...,p_{n-1}$ odd and $p_n$ even, revisited}
The goal of this subsection is to diagonlize the symmetrized linking matrix $\mathcal L + \mathcal L^\tau$
obtained in theorem \ref{linkingform1}. Specifically, we want to find a regular matrix $P$ of the same dimension 
as $\mathcal L + \mathcal L^\tau$ such that $P^\tau(\mathcal L + \mathcal L^\tau)P$ is a diagonal matrix. 
By way of shortcut of notation, we will write $\langle x, y\rangle$ to denote $\ell k (x,y) + \ell k(y,x)$. 

As the matrix $\mathcal L + \mathcal L^\tau$ from theorem \ref{linkingform1} consists of a number of matrix 
blocks of the form $\pm Y_m$ (see \eqref{X-n} for the definition of $Y_m$), we first take the time to apply the Gram-Schmidt process to the latter. 
We let $P_m$ denote the upper triangular  $m\times m$ matrix given by 
\begin{equation} \label{P-n}
P_m =  \left[ 
\begin{array}{rrrrrrr}
1 & -1 & -1 & -1 & ... & -1 & -1 \cr
0 & 2 & -1 & -1 & ... & -1 & -1 \cr
0 & 0 & 3 & -1 & ... & -1 & -1 \cr
\vdots & \vdots & \vdots & \vdots  & \vdots & \vdots & \vdots  \cr 
0 & 0 & 0 &  & ... & m-1 & -1  \cr
0 & 0 & 0 &  & ... & 0 & m \cr
\end{array}
\right]
\end{equation}
\begin{lemma} \label{aux1}
Consider the inner product space $(\langle \cdot , \cdot \rangle,\mathbb{Z}^m)$ where the inner product $\langle \cdot , \cdot \rangle$ with respect to the standard basis $\{\alpha _1, ..., \alpha _m\}$ of $\mathbb{Z}^m$ is given by 
$$\langle \alpha _i , \alpha _j\rangle = (i,j)\mbox{--th entry of the matrix } Y_m \mbox{ from \eqref{X-n}} $$
Then defining $a_i = i\alpha _i -\alpha _{i-1} -\alpha _{i-2} - ... - \alpha _1$ for each $i=1,...,m$ yields an orthogonal basis 
for $(\langle \cdot , \cdot \rangle,\mathbb{Z}^m)$ with $\langle a _i, a _i \rangle = i(i+1)$. 
Said differently,\footnote{Here and in the remainder of the article, 
we let Diag$\, (x_1,x_2,...,x_m)$ denote the $m\times m$ square matrix whose off-diagonal entries are zero and whose
diagonal entries are given by $x_1,...,x_m$.}

$$ P_m^\tau\,  Y_m \, P_m = \mbox{Diag}\, (1\cdot 2\, , \, 2\cdot 3 \, , \, 3\cdot 4 \, , ..., \, m\cdot (m+1))$$
\end{lemma}
\begin{proof}
This is a straightforward application of the Gram-Schmidt process. Let $a_i$ be as stated in the lemma and assume 
that $\{a_1,...,a_i\}$ is an orthogonal set for all $i<k\le  m$ with the stated squares $\langle a_i, a_i \rangle = i (i+1)$
(the case of $i=1$ being clearly true). We prove that the statement remains true if $i$ is chosen to be $k$. 
Note that 
$$\langle \alpha_k,a_i \rangle  = \langle \alpha _k,  i\alpha _i - \alpha _{i-1} - ... - \alpha _1 \rangle = 
i -1-1-...-1 = 1 $$
for any choice of $i<k$. Using the Gram-Schmidt process gives 
\begin{align} \nonumber 
a_k & = \alpha _k - \frac{\langle \alpha_k, a _{k-1}\rangle}{\langle a _{k-1}, a_{k-1}\rangle} a_{k-1} - ... - 
 \frac{\langle \alpha_k, a _{1}\rangle}{\langle a _{1}, a_{1}\rangle} a_{1} \cr
 & = \alpha _k - \frac{1}{(k-1)k}((k-1)\alpha _{k-1} - \alpha _{k-2}-...- \alpha _1) - ... - \alpha _1 \cr 
 & = \frac{1}{k}\left( k \alpha _k - \alpha _{k-1} - \alpha _{k-2} - ... - \alpha _1\right) \cr 
\end{align}
Proceeding as in remark \ref{gramgcd}, we let $a_k$ be equal to 
$$a_k =  k \alpha _k - \alpha _{k-1} - \alpha _{k-2} - ... - \alpha _1$$
which already showes that $\{a_1,...,a_k\}$ is orthogonal. To complete the proof of the lemma, we need to compute 
$\langle a_k, a_k \rangle$: 
\begin{align} \nonumber
\langle a_k, a_k \rangle & = \langle  k \alpha _k - \alpha _{k-1} - ... - \alpha _1,k \alpha _k - \alpha _{k-1} - ... - \alpha _1 \rangle \cr
& =  k^2 \langle \alpha _k , \alpha _k\rangle -2k \langle \alpha _k, \alpha _{k-1} + ... + \alpha _1 \rangle + 
\sum _{i=1}^{k-1} \langle \alpha _i , \alpha _i\rangle + \sum _{\tiny \begin{array}{c} i,j =1 \cr i\ne j \end{array}}^{k-1}
\langle \alpha _i , \alpha _j \rangle \cr 
& = 2k^2 -2k(k-1) + 2(k-1) + ((k-1)^2 - (k-1)) \cr 
& = k(k+1)
\end{align}
which is as claimed. 
\end{proof}

We proceed by defining vectors $a^i_k$ as 
$$ a^i_k = k\alpha ^i_k - \alpha ^i _{k-1} - ... - \alpha ^i_1 $$
for each $i=1,...,n-1$. 
Lemma \ref{aux1} then shows that for each such index $i$, the set $\{ a^i_1,...,a^i_{|p_i|-1}\}$ is an orthogonal set 
with respect to 
 $\lr = \mathcal L + \mathcal L^\tau$ and $\langle a^i_k, a^i_k \rangle = -Sign(p_i) \, k(k+1)$. Moreover, since $\langle \alpha ^i_k, \alpha ^j_m\rangle =0$
whenever $i\ne j$, we see that in fact the set 
\begin{equation} \label{b1prelim}
\mathcal B^\perp_{1,prelim} = \{a^1_1,...,a^1_{|p_1|-1},a^2_1,...,a^2_{|p_2|-1},...,a^{n-1}_1,...,a^{n-1}_{|p_{n-1}|-1}\}\end{equation}
is also an orthogonal set. 

We then turn to finding two additional vectors, which we shall label $X$ and $Y$, needed to complete $\mathcal B^\perp_{1,prelim}$
to an orthogonal basis for $\mathcal B_1^\perp$ for $H_1(\Sigma _1;\mathbb{Q})$. We find $X$ using again the Gram-Schmidt process.
%
%
\begin{lemma} \label{aux2}
Setting $X$ equal to 
$$X = |p_1 \cdot ... \cdot p_{n-1}| \,  \gamma - \sum _{i=1}^{n-1} \left(   \left( \prod _{k\ne i} |p_k| \right) \sum _{k=1}^{|p_i|-1} \alpha ^i_k \right)$$
makes the set $\mathcal B^\perp_{1,prelim}  \cup \{ X\}$ an orthogonal set. Moreover, the square of $X$ is 
$$ \langle X, X \rangle = - (p_1\cdot ... \cdot p_{n-1})^2 \left( \frac{1}{p_1} + ... + \frac{1}{p_{n-1}} \right) 
$$
\end{lemma}
\begin{proof}
An easy induction argument on $|p_i|$ shows that  
$$  \sum _{k=1}^{|p_i|-1}  \frac{\langle \gamma, a^i_k\rangle}{\langle a^i_k, a^i_k\rangle}  a^i_k = \frac{1}{|p_i|} 
(\alpha ^i_1 + \alpha ^i_2 + ... + \alpha ^i_{|p_i|-1} )$$
Letting $X$ be given by the Gram-Schmidt formula
$$ X = \gamma - \sum_{i=1}^{n-1} \sum _{k=1}^{|p_i|-1}  \frac{\langle \gamma, a^i_k\rangle}{\langle a^i_k, a^i_k\rangle} a^i_k$$
leads, in conjunction with the above formula, to 
$$ X = \gamma - \sum _{i=1}^{n-1} \frac{1}{|p_i|} \sum _{k=1}^{|p_i|-1} \alpha ^i_k $$
To keep coefficients integral (see remark \ref{gramgcd}) we multiply the right-hand side of the above by 
$|p_1\cdot ... \cdot p_{n-1}|$ and set $X$ instead equal to  
$$X = |p_1 \cdot ... \cdot p_{n-1}| \,  \gamma - \sum _{i=1}^{n-1} \left(   \left( \prod _{k\ne i} |p_k| \right) \sum _{k=1}^{|p_i|-1} \alpha ^i_k \right)$$
as in the statement of the lemma. Thus, $\mathcal B^\perp_{1,prelim}  \cup \{ X\}$ is indeed an orthogonal set.

We next compute $\langle X, X \rangle$:   
\begin{align} \nonumber
\langle X, X \rangle & = (p_1\cdot ... \cdot p_{n-1})^2 \langle \gamma,\gamma \rangle + \sum_{i=1}^{n-1}  \left[     \left( \prod _{k\ne i} |p_k| \right)^2 \left\langle  \sum _{k=1}^{|p_i|-1} \alpha ^i_k , \sum _{k=1}^{|p_i|-1} \alpha ^i_k   \right\rangle \right] - \cr
& -  2 \sum _{i=1}^{n-1} \frac{(p_1 \cdot ... \cdot p_{n-1})^2}{|p_i|} \left\langle \gamma, \sum _{k=1}^{|p_i|-1} \alpha ^i_k \right\rangle
\end{align}
Using the linking form $\mathcal L$ from theorem \ref{linkingform1}, it is easy to see that  (for example by induction 
on $|p_i|$) 
\begin{align} \nonumber
\left\langle  \sum _{k=1}^{|p_i|-1} \alpha ^i_k , \sum _{k=1}^{|p_i|-1} \alpha ^i_k   \right\rangle & = 
-Sign(p_i) \, |p_i| \, (|p_i| -1)  \cr
 \left\langle \gamma, \sum _{k=1}^{|p_i|-1} \alpha ^i_k \right\rangle & = - Sign(p_i) \, (|p_i|-1) 
\end{align}
which in turn shows that 
\begin{align} \nonumber
\sum_{i=1}^{n-1} &  \left[     \left( \prod _{k\ne i} |p_k| \right)^2 \left\langle  \sum _{k=1}^{|p_i|-1} \alpha ^i_k , \sum _{k=1}^{|p_i|-1} \alpha ^i_k   \right\rangle \right] -  2 \sum _{i=1}^{n-1} \frac{(p_1 \cdot ... \cdot p_{n-1})^2}{|p_i|} \left\langle \gamma, \sum _{k=1}^{|p_i|-1} \alpha ^i_k \right\rangle = \cr
= - & \sum_{i=1}^{n-1} Sign(p_i)  \frac{(p_1 \cdot ... \cdot p_{n-1})^2}{|p_i|} (|p_i|-1)  + 2 \sum_{i=1}^{n-1} Sign(p_i)  \frac{(p_1 \cdot ... \cdot p_{n-1})^2}{|p_i|} (|p_i|-1)  \cr
= & \sum_{i=1}^{n-1} Sign(p_i)  \frac{(p_1 \cdot ... \cdot p_{n-1})^2}{|p_i|} (|p_i|-1) \cr
= & |p_1\cdot ... \cdot p_{n-1}| \,  \sum_{i=1}^{n-1} Sign(p_i) \left( \prod_{k\ne i} |p_k|\right)  (|p_i|-1)
\end{align}
Finally, recalling (see theorem \ref{linkingform1}) that $\langle \gamma , \gamma \rangle = - (Sign(p_1) + ... + Sign(p_{n-1}))$, we are able to assemble all the pieces to compute $\langle X, X \rangle$:
\begin{align} \nonumber
\frac{1}{|p_1\cdot ... \cdot p_{n-1}|}\,  \langle X, X \rangle & = -|p_1\cdot ... \cdot p_{n-1}|(Sign(p_1) + ... + Sign(p_{n-1})) + \cr
 & \quad \quad \quad \quad +  \sum_{i=1}^{n-1} Sign(p_i) \left( \prod_{k\ne i} |p_k|\right)  (|p_i|-1) \cr
& = - \sum _{i=1}^{n-1}  Sign(p_i) \prod _{k\ne i } |p_k| 
\end{align}
and so 
$$ \langle X, X \rangle =  - (p_1\cdot ... \cdot p_{n-1})^2 \cdot \left(\frac{1}{p_1} + ... + \frac{1}{p_{n-1}}  \right)$$
as claimed in the statement of the lemma.
\end{proof}
In the final step, we would like to find a vector $Y\in H_1(\Sigma _1;\mathbb{Q})$ such that $\mathcal B_{1,prelim}^\perp \cup \{X,Y\}$ is an orthogonal basis. While $\langle a^i_k, a^i_k \rangle \ne 0$ for any choice of $i,k$, and thus the Gram-Schmidt process worked well for finding $X$, it is possible, and 
it does happen, that $\langle X, X \rangle =0$. This of course obstructs us from finding $Y$ by means of the Gram-Schmidt process, calling instead for an application of theorem \ref{squarezero}. 
We proceed by treating the two cases $\langle X, X, \rangle \ne 0$ and $\langle X, X, \rangle = 0$ separately. 
\begin{lemma} \label{aux3}
Let $X\in H_1(\Sigma _1 ; \mathbb{Q})$ be as defined in lemma \ref{aux2} and assume that $\langle X, X \rangle \ne 0$. Define $Y\in H_1(\Sigma _1 ; \mathbb{Q})$ as 
\begin{align} \nonumber
Y & = |p_1\cdot ... \cdot p_{n-1}| \,\left( \frac{1}{p_1} + ... + \frac{1}{p_{n-1}} \right) \delta +X \cr
\end{align}
Then $\mathcal B_1^\perp = \mathcal B^\perp _{1,prelim} \cup \{X,Y\}$ is an orthogonal basis and 
\begin{align} \nonumber
\langle Y, Y \rangle = 
 \left( \prod _{i=1}^{n-1} p_1\cdot ... \cdot \hat p_i \cdot ... \cdot p_{n-1} \right) \cdot  \left( \prod _{i=1}^n p_1\cdot ... \cdot \hat p_i \cdot ... \cdot p_n \right)
\end{align}
\end{lemma}
\begin{proof}
Our assumption $\langle X, X \rangle \ne 0$ allows us to use the Gram-Schmidt process to find $Y$ as 
\begin{align} \nonumber
Y= \delta - \frac{\langle \delta, X\rangle}{\langle X, X\rangle} X - \sum_{i=1}^{n-1} \sum _{k=1}^{|p_i|-1}  \frac{\langle \delta, a^i_k\rangle}{\langle a^i_k, a^i_k\rangle} a^i_k 
\end{align}
Since $\langle \delta, \alpha ^i_k\rangle =0$ for all $i,k$ it follows that $\langle \delta, a ^i_k\rangle =0$ also, reducing the 
above formula to 
$$Y= \delta - \frac{\langle \delta, X\rangle}{\langle X, X\rangle} X$$
With $\langle X, X \rangle$ already computed in lemma \ref{aux2}, the same lemma (using also the result of theorem \ref{linkingform1}) implies that  
$$ \langle \delta, X \rangle = |p_1\cdot ... \cdot p_{n-1} |$$
showing that 
\begin{align} \nonumber
Y & = \delta + \frac{1}{|p_1\cdot ... \cdot p_{n-1}| \, \left( \frac{1}{p_1} + ... + \frac{1}{p_{n-1}} \right) } \, X 
\end{align} 
To keep our coefficients integral (see remark \ref{gramgcd}) we instead set 
$$ 
Y = |p_1\cdot ... \cdot p_{n-1}| \, \left( \frac{1}{p_1} + ... + \frac{1}{p_{n-1}} \right) \delta +X 
$$
showing that $\mathcal B_1^\perp = \mathcal B^\perp _{1,prelim} \cup \{X,Y\}$ is an orthogonal basis for $(\lr,H_1(\Sigma_1;\mathbb{Q}))$. It remains to calculate $\langle Y, Y \rangle$:
\begin{align} \nonumber
\langle Y, Y \rangle & = (p_1\cdot ... \cdot p_{n-1})^2  \, \left( \frac{1}{p_1} + ... + \frac{1}{p_{n-1}} \right) ^2 \langle \delta, \delta \rangle + \cr 
& \quad \quad \quad \quad +2  |p_1\cdot ... \cdot p_{n-1}| \, \left( \frac{1}{p_1} + ... + \frac{1}{p_{n-1}} \right) \langle \delta, X\rangle + \langle X, X \rangle \cr
& = (p_1\cdot ... \cdot p_{n-1})^2  \, \left( \frac{1}{p_1} + ... + \frac{1}{p_{n-1}} \right) ^2 p_n + \cr
& \quad \quad \quad \quad  + 2  (p_1\cdot ... \cdot p_{n-1})^2 \, \left( \frac{1}{p_1} + ... + \frac{1}{p_{n-1}} \right) 
- ( p_1\cdot ... \cdot p_{n-1}) ^2 \left( \frac{1}{p_1} + ... + \frac{1}{p_{n-1}} \right) \cr 
& = (p_1\cdot ... \cdot p_{n-1})^2  \, \left( \frac{1}{p_1} + ... + \frac{1}{p_{n-1}} \right) \left[ \left( \frac{1}{p_1} + ... + \frac{1}{p_{n-1}} \right)p_n + 1 \right]
\end{align}
\end{proof}
\begin{lemma} \label{aux4}
Let $X\in H_1(\Sigma _1 ; \mathbb{Q})$ be as defined in lemma \ref{aux2} and assume that $\langle X, X \rangle = 0$. Then, in the Witt ring $W(\mathbb{Q})$, the equality 
$$ ( \langle \cdot , \cdot \rangle, H_1(\Sigma _1;\mathbb{Q})) = (\langle \cdot , \cdot \rangle |_{V\times V},V)$$
holds where $V=Span \, \mathcal B^\perp _{1,prelim}$ (where the latter is as defined in \eqref{b1prelim}). 
\end{lemma}
\begin{proof}
This is a direct consequence of theorem \ref{squarezero} and can also be verified directly. 
Namely, observe that the format of $\mathcal L + \mathcal L^\tau$ as 
calculated in theorem \ref{linkingform1} shows that
$$ ( \langle \cdot , \cdot \rangle, H_1(\Sigma _1;\mathbb{Q})) = 
 (\langle \cdot , \cdot \rangle |_{V\times V},V) \oplus  (\langle \cdot , \cdot \rangle |_{W\times W},W)$$
where $W = Span \, \{ X, \delta\}$. But since $\langle \cdot , \cdot \rangle |_{W\times W}$ is represented 
by the matrix 
$$ \langle \cdot , \cdot \rangle |_{W\times W} = \left[
\begin{array}{cc}
0 & |p_1\cdot ... \cdot p_{n-1} |  \cr 
|p_1\cdot ... \cdot p_{n-1} |  & p_n
\end{array}
\right]$$
with respect to the basis $\{X,\delta\}$, we see that $\langle \cdot , \cdot \rangle |_{W\times W}$ is metabolic and 
thus equivalent to zero in $W(\mathbb{Q})$. 
\end{proof}
We summarize our findings in the next theorem:
\begin{theorem}  \label{diagmat1}
Let $P$ be the upper triangular matrix 
$$P = \left[ 
\begin{array}{ccc|c|ccc|c|c}
   &   &       &  \dots   &    &       &     &  \hat p_1 p_2  ... p_{n-1} & \hat p_1 p_2  ... p_{n-1}  \cr 
   & P_{|p_1|-1} & &     \dots      &    &    &     &  \vdots & \vdots  \cr 
   &  & &      \dots  &    &   &         &  \hat p_1 p_2  ... p_{n-1} & \hat p_1 p_2  ... p_{n-1}  \cr \hline 
  \vdots  &  \vdots   &  \vdots   &  \ddots  &    &   &         &  \vdots  &   \vdots \cr \hline
   &  & &      &    &   &         & p_1 p_2  ... \hat p_{n-1}  & p_1 p_2  ... \hat p_{n-1}   \cr 
   &  & &      &    & P_{|p_{n-1}|-1}  &         & \vdots  & \vdots   \cr 
   &  & &      &    &   &         & p_1 p_2  ... \hat p_{n-1}  & p_1 p_2  ... \hat p_{n-1}   \cr  \hline 
  0 & \dots &0     & \dots  &  0 & \dots  & 0        &  |p_1\cdot ... \cdot p_{n-1}| & |p_1\cdot ... \cdot p_{n-1}|   \cr \hline 
   0 & \dots  & 0     & \dots  &   0 &  \dots  &    0     & 0 &  |p_1\cdot ... \cdot p_{n-1}| \cdot \left( \sum _{i=1}^{n-1} \frac{1}{p_i} \right) \cr 
\end{array} 
\right]
$$
where $P_m$ is as defined in \eqref{P-n} and let $\mathcal L$ be as computed in theorem \ref{linkingform1}. Then, if $\langle X, X \rangle \ne 0$, one gets 
\begin{align} \nonumber
P^\tau  (\mathcal L & + \mathcal L^\tau) P   = \cr
 \oplus & \left( \bigoplus _{i=1}^{n-1} Diag (-Sign(p_i) \cdot  1\cdot 2, -Sign(p_i)  \cdot 2\cdot 3, ..., -Sign(p_i) \cdot (|p_i|-1)\cdot |p_i|) \right)  \oplus \cr
\oplus & Diag (\langle X, X \rangle, \langle Y, Y \rangle   ) 
\end{align}
If $\langle X, X \rangle = 0$, let $Q$ be the matrix obtained from $P$ by setting its last column and row 
equal to zero, safe the diagonal entry which should be set equal to 1. 
Then  
\begin{align} \nonumber
Q^\tau  (\mathcal L & + \mathcal L^\tau) Q   = \cr
 \oplus & \left( \bigoplus _{i=1}^{n-1} Diag (-Sign(p_i) \cdot  1\cdot 2, -Sign(p_i)  \cdot 2\cdot 3, ..., -Sign(p_i) \cdot (|p_i|-1)\cdot |p_i|) \right) \oplus \cr
 \oplus &  \left[\begin{array}{cc} 0 & |p_1\cdot ... \cdot p_{n-1}| \cr |p_1\cdot ... \cdot p_{n-1}| & p_n \end{array}\right]
\end{align}
Recall that $\langle X, X \rangle $ and $\langle Y, Y \rangle$ have been computed in lemmas \ref{aux2} and \ref{aux3}.
\end{theorem}
Before continuing on, we take a moment to express the quantities  $\langle X, X \rangle$ and $\langle Y, Y \rangle$ in more familiar terms involving determinants of knots/links. 
%
%
%
\begin{lemma} \label{aux5}
Assume that $n,p_1,...,p_{n-1}$ are odd integers with $n\ge 3$ and that $p_n\ne 0$ is an even integer. 
Consider the pretzel knot $P(p_1,...,p_n)$ and the pretzel link (of 2 components) $P(p_1,..,p_{n-1})$. 
Then 
\begin{align} \nonumber 
\det P(p_1,...,p_n) & =  \prod _{i=1}^n p_1 \cdot ... \cdot \hat p _i \cdot ... \cdot p_n  \cr
\det P(p_1,...,p_{n-1}) & = \prod _{i=1}^{n-1} p_1 \cdot ... \cdot \hat p _i \cdot ... \cdot p_{n-1} 
\end{align}
In particular, we can re-write $\langle X, X \rangle$ and $\langle Y,Y\rangle$ as 
\begin{align} \nonumber
\langle X,X\rangle & = -(p_1 \cdot ... \cdot p_{n-1})  \cdot \det P(p_1,...,p_{n-1}) \cr 
\langle Y,Y\rangle & = \det  P(p_1,...,p_n) \cdot \det P(p_1,...,p_{n-1})
\end{align}
\end{lemma}
\begin{proof}
We shall calculate $\det P(p_1,...,p_n)$ by relying on the formula $\det P(p_1,...,p_n)= \det (\mathcal L -  \mathcal L^\tau)$ with $\mathcal L$ as in theorem \ref{linkingform1}. As we shall see, the result of this computation 
agrees with the formula claimed by the lemma only up to sign. We allow ourselves the liberty of choosing the 
sign of the determinant somewhat arbitrarily.

If $\langle X, X \rangle \ne 0$, we simply apply the determinant to the relation $A=P^\tau (\mathcal L + \mathcal L^\tau)P$ from theorem \ref{diagmat1} (where we let $A$ denote the first diagonal matrix from that theorem): 
\begin{align} \nonumber
\det (\mathcal L + \mathcal L^\tau) & =  \det A/ ( \det P )^2  \cr
& = - \frac{\left( \prod_{i=1}^{n-1}\frac{(|p_i|!)^2}{|p_i|}\right)(p_1\cdot ... \cdot p_{n-1})^4 \left( \frac{1}{p_1}+ ... + \frac{1}{p_{n-1}} \right)^2 \cdot \left[\left( \frac{1}{p_1}+ ... + \frac{1}{p_{n-1}} \right)p_n +1     \right] }{\left( \prod _{i=1}^{n-1} \frac{|p_i|!}{|p_i|} \right)^2 \cdot (p_1\cdot ... \cdot p_{n-1})^4 \cdot 
\left( \frac{1}{p_1}+ ... + \frac{1}{p_{n-1}} \right)^2} \cr
& = -|p_1\cdot ... \cdot p_{n-1}|\cdot \left[ \left( \frac{1}{p_1}+ ... + \frac{1}{p_{n-1}} \right)p_n +1 \right] \cr
& = - Sign (p_1\cdot ... \cdot p_{n-1}) \cdot \left( \prod _{i=1}^n p_1\cdot ... \cdot \hat p_i \cdot ... \cdot p_n  \right)
\end{align}
If $\langle X, X \rangle = 0$ a similar argument applies. Namely, applying the determinant to 
the equation $Q^\tau (\mathcal L + \mathcal L^\tau)Q$ from theorem \ref{linkingform1}, yields   the 
desired result, the details are left as an easy exercise. 
\vskip3mm
The computation of $\det P(p_1,...,p_{n-1})$ 
follows along the same lines with only minor modification. We focus on these differences rather than repeating
the entire calculation. 

The reader should first note that the Seifert surface $\Sigma_1$ for $P(p_1,...,p_n)$ displayed in figure \ref{pic4}, 
becomes a Seifert surface for $P(p_1,...,p_{n-1})$ after removing the unique band with $p_n$ half twists. 
We shall call the resulting surface $\Sigma_1 '$. Its linking form $\mathcal L'$ differs from $\mathcal L$ only in the 
last row and column (which are removed from $\mathcal L$ to obtain $\mathcal L'$).  In particular, the computation 
of $\det (\mathcal L' + \mathcal L'^\tau)$ is identical to that of $\det (\mathcal L + \mathcal L^\tau)$ safe 
the contribution of $Y$ to the latter. Thus, 
\begin{align} \nonumber
\det (\mathcal L '+ \mathcal L'^\tau) & = \det (\mathcal L + \mathcal L^\tau) \,  \frac{(\mbox{Coefficient of $\delta$ in $Y$})^2}{\langle Y, Y \rangle }  \cr
&    =  - Sign (p_1\cdot ... \cdot p_{n-1}) \left( \prod _{i=1}^n p_1\cdot ... \cdot \hat p_i \cdot ... \cdot p_n \right) \cdot \cr
& \quad \quad \quad \quad \quad \quad \quad \quad  \cdot 
\frac{\left( \prod _{i=1}^{n-1} p_1 \cdot ... \cdot \hat p_i \cdot ... \cdot p_{n-1} \right)^2}{\left( \prod _{i=1}^{n-1} p_1\cdot ... \cdot \hat p_i \cdot ... \cdot p_{n-1} \right) \cdot  \left( \prod _{i=1}^n p_1\cdot ... \cdot \hat p_i \cdot ... \cdot p_n \right)} \cr
& =  - Sign (p_1\cdot ... \cdot p_{n-1}) \left( \prod _{i=1}^{n-1} p_1\cdot ... \cdot \hat p_i \cdot ... \cdot p_{n-1}\right)
\end{align} 
This formula applies in both the cases when $\langle X, X\rangle=0$ and $\langle X, X\rangle\ne 0$. With this 
observation, the proof of the lemma is complete. 
\end{proof}
%
%
%
\subsection{The case of $n, p_n$ even and $p_1,...,p_{n-1}$ odd, revisited}
In this section we turn to diagonalizing the symmetrized linking form $\mathcal L + \mathcal L^\tau$ 
with $\mathcal L$ this time being as computed in theorem \ref{linkingform2}. The work has largely been 
done in the previous section and we focus our attention only on the minor differences. 
%
%
%
\begin{lemma} \label{detcompu2}
Let $n$ and $p_n$ be even integers with $n\ge 3$ and $p_n\ne 0$ and let $p_1,...,p_{n-1}$ be odd integers. 
Let $\mathcal L$ be the linking matrix associated to the Seifert surface $\Sigma _3$ of $P(p_1,...,p_n)$ 
and the basis $\mathcal B_3$ of $H_1(\Sigma _3;\mathbb{Q})$ as defined in figure \ref{pic7}. Then the 
determinant of $P(p_1,...,p_n)$ is 
$$ \det P(p_1,...,p_n) =  \prod _{i=1}^n p_1 \cdot ... \cdot \hat p_i \cdot ... \cdot p_n  
$$
\end{lemma}
\begin{proof}
Recall that $\det P(p_1,...,p_n) = \det (\mathcal L + \mathcal L^\tau)$ but that we allow ourselves the freedom 
of choosing the sign of the determinant. 

The determinant of $ \mathcal L + \mathcal L^\tau $ is computed in analogy to the computation of lemma \ref{aux5}. 
Specifically, let $P'$ be the matrix obtained from the matrix $P$ from theorem \ref{linkingform1} by deleting its last row 
and column, and let $A'$ be the diagonal matrix $P^\tau (\mathcal L + \mathcal L')P$  from the same theorem, again with its last row and 
column deleted. Then $(P')^\tau  (\mathcal L + \mathcal L^\tau) P' = A'$ so that 
\begin{align} \nonumber 
\det \mathcal L + \mathcal L^\tau & = \frac{\det A'}{(\det P')^2} \cr 
& = Sign(p_n)  \frac{\left( \prod_{i=1}^{n}\frac{(|p_i|!)^2}{|p_i|}\right)(p_1\cdot ... \cdot p_{n})^2 \left( \frac{1}{p_1}+ ... + \frac{1}{p_{n}} \right)    }{\left( \prod _{i=1}^{n} \frac{|p_i|!}{|p_i|} \right)^2 \cdot (p_1\cdot ... \cdot p_{n})^2 } \cr
& = Sign(p_n) |p_1\cdot ... \cdot p_{n}|\cdot \left( \frac{1}{p_1}+ ... + \frac{1}{p_{n}} \right)    \cr
& = Sign (p_1\cdot ... \cdot p_{n-1}) \cdot \left( \prod _{i=1}^n p_1\cdot ... \cdot \hat p_i \cdot ... \cdot p_n  \right) 
\end{align}
as needed, up to sign. 
\end{proof}
We have thus proved the following theorem:
\begin{theorem} \label{diagmat2}
Let $n$ and $p_n$ be even integers with $n\ge 3$ and $p_n\ne 0$ and let $p_1,...,p_{n-1}$ be odd integers. 
Let $P$ be the matrix 
$$P = \left[ 
\begin{array}{ccc|c|ccc|c}
   &   &       &  \dots   &    &       &     &  \hat p_1 p_2  ... p_{n}   \cr 
   & P_{|p_1|-1} & &     \dots      &    &    &     &  \vdots   \cr 
   &  & &      \dots  &    &   &         &  \hat p_1 p_2  ... p_{n}  \cr \hline 
  \vdots  &  \vdots   &  \vdots   &  \ddots  &    &   &         &  \vdots   \cr \hline
   &  & &      &    &   &         & p_1 p_2  ... \hat p_{n}     \cr 
   &  & &      &    & P_{|p_{n}|-1}  &         & \vdots     \cr 
   &  & &      &    &   &         & p_1 p_2  ... \hat p_{n}     \cr  \hline 
  0 & \dots &0     & \dots  &  0 & \dots  & 0        &  |p_1\cdot ... \cdot p_{n}|  
\end{array} 
\right]
$$
and let $\mathcal L$ be as computed in theorem \ref{linkingform2}. Then 
\begin{align} \nonumber
P^\tau (\mathcal L + \mathcal L^\tau) P & = \left(  \bigoplus _{i=1}^{n}  Diag  ( -Sign(p_i) \cdot  1\cdot 2  ,  ... ,    -Sign(p_i) \cdot (|p_i|-1)\cdot |p_i|  ) \right) \oplus \cr
& \quad\quad \quad\quad\quad \quad \quad \quad  \oplus  Diag(   -(p_1\cdot ... \cdot p_n) \cdot  \det P(p_1,...,p_n) ) 
\end{align}
The determinant $\det P(p_1,...,p_n) $  has been computed in lemma \ref{detcompu2}. 
\end{theorem}
%
%
%
\subsection{The case of $n, p_1,...,p_n$ all odd, revisited}
The goal of this section is to diagonalized the symmetrized linking matrix $\mathcal L + \mathcal L^\tau$ 
from theorem \ref{linkingform3}. Here too we would like to utilize the Gram-Schmidt process inasmuch as possible. 
Recall that the basis $\mathcal B_3$ for $H_1(\Sigma _3;\mathbb{Q})$ is $\mathcal B_3=\{\alpha _1, ...,\alpha _{n-1}\}$  with $\alpha_i$ as in figure \ref{pic6}. We wish to create an orthogonal basis 
$\mathcal B_3^\perp = \{ a_1,...,a_{n-1}\}$ by means of the formalism from theorem \ref{gramschmidttheorem} (see also remark \ref{gramgcd}). Towards that goal, we prove a simple lemma after reminding the reader of some notation which was mentioned in the introduction. 

For an integer $i\ge 1$, let $\sigma_i(t_1,...,t_m)$ be the $i$-th symmetric polynomial in the variables 
$t_1,...,t_m$. For example, $\sigma _1(t_1,...,t_n) = t_1+...+t_m$ and 
$\sigma _2(t_1,...,t_m) = t_1t_2 + t_1t_3+ ... + t_{m-1}t_m$ and so on.  By convention, we define the 0-th symmetric polynomial to be $\sigma_0(t_1,...,t_m)=1$. We shall write $\sigma _i$ for $\sigma_i(p_1,...,p_{i+1})$. 
\begin{lemma} \label{auxxx1}
Set $a_1=\alpha _1$ and $a_{i+1} = \sigma_{i}\, \alpha _{i+1} + p_{i+1} a_i \in H_1(\Sigma _3;\mathbb{Q})$ and let 
 $\mathcal B_3^\perp = \{a_1,...,a_{n-1}\}$. Then $\mathcal B_3^\perp$ is an orthogonal set  and   
$$ \langle a_i , a_i \rangle  = \sigma_{i-1} \cdot \sigma _{i} $$ 
\end{lemma}
Before proving this statement, we would like to point out that lemma \ref{auxxx1} does not claim, indeed this would be false in certain cases, that $\mathcal B_3^\perp$ is a basis for $H_1(\Sigma _3;\mathbb{Q})$. Some elements of 
$\mathcal B_3^\perp$ may be zero. 
\begin{proof} 
We proof lemma \ref{auxxx1} by induction on $i$, the cases of $i=1,2$ are easily seen to hold. Proceeding 
to the step of the induction, we consider the vector $a_{i+1}$. Pick first an index $j$ with $j<i$, then we get 
\begin{align} \nonumber
\langle a_{i+1} , a_j \rangle & = \langle  \sigma_{i} \, \alpha _{i+1}, a_j \rangle  + \langle p_{i+1} a_i,
a_j \rangle = 0 
\end{align} 
since in this case $\langle \alpha _{i+1},a_j\rangle =0$ (as follows from inspection of the linking matrix $\mathcal L$ from theorem \ref{linkingform3}). On the other hand, 
\begin{align} \nonumber
\langle a_{i+1} , a_i \rangle & = \langle  \sigma_{i}\, \alpha _{i}, a_i \rangle  + \langle p_{i+1} a_i,a_i \rangle \cr
& =  -\sigma_{i} \cdot \sigma_{i-1} \cdot p_{i+1} + p_{i+1} \cdot \sigma_{i-1} \cdot \sigma _{i}\cr
& = 0 
\end{align} 
To finish the induction argument, we next determine $\langle a_{i+1}, a_{i+1} \rangle$: 
\begin{align} \nonumber
\langle a_{i+1}, a_{i+1} \rangle & = \langle  \sigma_{i}\, \alpha _{i+1} + p_{i+1} a_i, \sigma_{i}\, \alpha _{i+1} + p_{i+1} a_i \rangle \cr
& = (\sigma_{i})^2 (p_{i+1}+p_{i+2}) +  p_{i+1}^2 \sigma_{i-1} \, \sigma _{i} - 2p^2_{i+1}  \sigma_{i} \, \sigma_{i-1} \cr
& = \sigma_{i}\, [ \sigma_{i}\, (p_{i+1}+p_{i+2}) - p_{i+1}^2\sigma_{i-1} ] \cr
& = \sigma_{i}\, \sigma_{i+1} 
\end{align} 
In the second to last line, we relied on the easy to verify identities
\begin{align} \nonumber
\sigma _j & = p_{j+1}\sigma _{j-1} + \sigma _j(p_1,...,p_j) \cr
\sigma _{j+1} & = p_{j+1}p_{j+2}\sigma _{j-1} + (p_{j+1}+p_{j+2})\sigma _j(p_1,...,p_j) 
\end{align}
\end{proof}

As the proof of lemma \ref{auxxx1} shows, the Gram-Schmidt algorithm breaks down whenever 
$\sigma_i$ vanishes for some $i\ge 1$. 
\begin{theorem} \label{diagmat3}
Let $n,p_1,...,p_n$ be odd integers with $n\ge 3$. Let $\mathcal L$ be the linking matrix for the pretzel 
knot $P(p_1,...,p_n)$ as described in theorem \ref{linkingform3}. Let $P$ be the upper triangular matrix 
$$P = \left[ 
\begin{array}{cccccccc}
\sigma_0 & b_{1,2}  &  b_{1,3}     &  \dots   &   b_{1,n-1}   \cr 
0 & \sigma_1  & b_{2,3}   &  \dots   &   b_{2,n-1}    \cr   
0 & 0  & \sigma_2      &  \dots   &   b_{3,n-1}    \cr 
\vdots & \vdots & \vdots  & ... & \vdots \cr
0 & 0 & 0 & ... & \sigma _{n-1} 
\end{array} 
\right] \quad \mbox{ with } \quad b_{k,i} = \sigma _{k-1} \cdot \prod _{j=k+1}^i p_j 
$$
Then 
$$ P^\tau (\mathcal L + \mathcal L^\tau) P = Diag ( \sigma _0 \cdot \sigma _1\, , \, \sigma _1\cdot \sigma_2\, ,... ,\,
\sigma_{n-2}\cdot \sigma _{n-1} )$$
The rational Witt class of $\mathcal L + \mathcal L^\tau$ is given by
$$ \varphi (P(p_1,...,p_n)) = \bigoplus _{i=1}^{n-1} \langle \sigma _{i-1} \cdot \sigma _i \rangle $$ 
\end{theorem}
\begin{proof}
The claim about the form of $P^\tau (\mathcal L + \mathcal L^\tau) P$ follows directly from 
lemma \ref{auxxx1}. The fact that the integers $b_{k,i}$ take the form described, can be proved by 
induction on $i$ by using the formulae (the first two lines being the definitions of $b_{*,*}$ as change of basis parameters, the third line being from lemma \ref{auxxx1})
\begin{align} \nonumber
a_{i+1} & =\sigma_i \alpha _{i+1} + b_{i,i+1}\alpha _i + b_{i-1,i+1}\alpha _{i-1} +  ... b_{1,i+1}\alpha _1 \cr
a_{i} & =\sigma_{i-1} \alpha _{i} + b_{i-1,i}\alpha _{i-1} +b_{i-2,i}\alpha _{i-2} +  ... b_{1,i}\alpha _1 \cr
a_{i+1} &  =  \sigma_i \alpha _{i+1} + p_{i+1} a_i
\end{align}

The claim of the theorem about Witt classes follows immediately from lemma \ref{auxxx1} in the case when none of the 
numbers $\sigma_i$ vanish since in that case the set $\mathcal B_3^\perp$ from the said lemma 
is actually a basis for $H_1(\Sigma _3;\mathbb{Q})$. We thus need to address the case when some of the 
$\sigma_i$ equal zero.
We shall prove the theorem by induction on $n$. 

When $n=3$ the symmetrized linking matrix $\mathcal L + \mathcal L ^\tau$ looks like 
$$\mathcal L + \mathcal L^\tau  = \left[ 
\begin{array}{cc}
p_1+p_2 & -p_2 \cr
-p_2 & p_2+p_3 
\end{array}
\right]
$$ 

If $\sigma _1 = p_1+p_2$ vanishes then $\mathcal L +\mathcal L^\tau$ is metabolic and thus zero 
in $W(\mathbb{Q})$. Conversely, if $p_1+p_2=0$ then $\langle \sigma _0 \sigma_1 \rangle \oplus 
\langle \sigma _1 \sigma _2\rangle  = 0 \in W(\mathbb{Q})$. If on the other hand $\sigma _2=p_1p_2+p_1p_3+p_2p_3$ vanishes (but $p_1+p_2$ does not), then the matrix representing $\langle \cdot , \cdot \rangle$ with respect to the basis $\{a_1, a_2\}$  is 
$$ \left[ 
\begin{array}{cc}
p_1+p_2 & 0 \cr
0 & 0 
\end{array}
\right] = Diag(p_1+p_2 , 0 ) 
$$ 
so that in this case $\mathcal L + \mathcal L^\tau$ equals $\langle p_1+p_2\rangle$ in $W(\mathbb Q)$. But, with 
the same vanishing assumption, we also get  $\langle \sigma _0 \sigma_1 \rangle \oplus 
\langle \sigma _1 \sigma _2\rangle  = \langle \sigma _0 \sigma_1 \rangle \oplus 
\langle 0 \rangle  = \langle p_1+p_2\rangle \in W(\mathbb{Q})$. This proves the theorem for the case of $n=3$. 

 To address the step of the induction, let $i$ be the smallest index for which $\sigma _i(p_1,...,p_{i+1})$ vanishes and consider the basis 
$\{a_1,a_2,...,a_i,\alpha _{i+1},...,\alpha _{n-1}\}$. Note that then $\langle a_i, a_i \rangle =0$. With respect to 
this basis, the intersection form $\langle \cdot, \cdot \rangle$ is represented by the matrix 
$$Diag(\sigma_0 \sigma _1, \sigma_1 \sigma _2, ..., \sigma_{i-2} \sigma _{i-1}) \oplus   \left[
\begin{array}{ccccc}
0  &  -p_{i+1}\sigma _{i-1} &  0 & 0 & \dots  \cr
-p_{i+1}\sigma _{i-1} & p_{i+1}+p_{i+2} & -p_{i+2} &  0 & \dots  \cr
0 & -p_{i+2}  &  p_{i+2}+p_{i+3} &  -p_{i+3}  & \dots  \cr
0 & 0 & -p_{i+3} & p_{i+3}+p_{i+4} & \dots   \cr
\vdots & \vdots & \vdots & \vdots & \ddots   \cr
\end{array}
\right]
$$
Consider the second of these two matrix summands. Add the first row multiplied by $-p_{i+2}/(p_{i+1}\sigma _{i-1})$
to the third row and likewise add the first column multiplied by $-p_{i+2}/(p_{i+1}\sigma _{i-1})$ to the third 
column (this simply corresponds to another change of basis). Thus we see that $\mathcal L + \mathcal L^\tau$ is 
represented by the matrix 
\begin{align}\nonumber
Diag(\sigma_0 \sigma _1, \sigma_1 \sigma _2, ..., \sigma_{i-2} \sigma _{i-1})  & \oplus   
\left[
\begin{array}{cc}
0  &  -p_{i+1}\sigma _{i-1}  \cr
-p_{i+1}\sigma _{i-1} & p_{i+1}+p_{i+2} 
\end{array}
\right]
\oplus \cr
& \oplus  \left[
\begin{array}{ccc}
  p_{i+2}+p_{i+3} &  -p_{i+3}  & \dots  \cr
 -p_{i+3} & p_{i+3}+p_{i+4} & \dots   \cr
 \vdots & \vdots & \ddots   \cr
\end{array}
\right]
\end{align}
The second summand is metabolic and therefore zero in $W(\mathbb{Q})$. On the third summand however
can apply the induction hypothesis and we conclude that  
\begin{align} \nonumber
\varphi(P(p_1,...,p_n)) & = Diag(\sigma_0 \sigma _1, \sigma_1 \sigma _2, ..., \sigma_{i-2} \sigma _{i-1}) \oplus \cr
&\oplus \left(  \bigoplus _{j=i+2}^{n-1} \langle \sigma _{j-i-2}(p_{i+2},...,p_{j})\sigma _{j-i-1}(p_{i+2},...,p_{j+1})\rangle
\right)   
\end{align}  
It remains to compare this to the result claimed by the theorem. For this purpose we observe that for $k\ge i$, the 
equality 
$$ \sigma _k(p_1,...,p_{k+1}) = \sigma _i(p_1,...,p_{i+1}) \, \sigma _{k-i}(p_{i+2},...,p_{k+1}) + 
\sigma _{i+1}(p_1,...,p_{i+1}) \, \sigma _{k-i-1}(p_{i+2},...,p_{k+1}) $$
holds. Thus in the event when $\sigma _i = 0$ we get that 
$$ \sigma _k(p_1,...,p_{k+1}) =\sigma _{i+1}(p_1,...,p_{i+1}) \, \sigma _{k-i-1}(p_{i+2},...,p_{k+1}) $$
Therefore, for $k\ge i+2$ we also get 
\begin{align} \nonumber
\langle \sigma _{k-1} \sigma _{k} \rangle & = \langle (\sigma _{i+1}(p_1,...,p_{i+1}))^2  \, \sigma _{k-i-2}(p_{i+2},...,p_{k}) \, \sigma _{k-i-1}(p_{i+2},...,p_{k+1})  \rangle \cr
& = \langle \sigma _{k-i-2}(p_{i+2},...,p_{k+1}) \, \sigma _{k-i-1}(p_{i+2},...,p_{k+1})  \rangle
\end{align}
while of course for $k=i,i+1$ we get $\langle \sigma _{k-1}\sigma _k\rangle = 0 \in W(\mathbb{Q})$. This completes
the proof of the induction step and thus of the theorem.
\end{proof}
\begin{lemma} 
Assume that $n,p_1,...,p_n$ are all odd with $n\ge 3$. Then the determinant of $P(p_1,...,p_n)$ is given 
by 
$$\det P(p_1,...,p_n) = \prod _{i=1}^{n} p_1 \cdot ... \cdot \hat p_i \cdot ... \cdot p_n  = \sigma _{n-1}$$ 
\end{lemma}
\begin{proof}
Let $\mathcal L$ be the linking matrix (from theorem \ref{linkingform3}) for $P(p_1,...,p_n)$ associated to the Seifert surface $\Sigma _3$ and the choice of basis $\mathcal B_3$ as in figure \ref{pic6}. 
%
%
%
%
%

To compute $\det (\mathcal L + \mathcal L^\tau)$ we proceed by induction on $n$. When $n=3$, 
the explicit form of $\mathcal L + \mathcal L^\tau$ from theorem \ref{linkingform3} shows that 
$\det (\mathcal L + \mathcal L^\tau) = p_1p_2+p_1p_3+p_2p_3$ as claimed by the lemma. When $n>3$, 
let $Y_n=Y_n(p_1,...,p_n)$ denote the matrix $\mathcal L + \mathcal L^\tau$ from theorem \ref{linkingform3}
but temporarily allowing $n$ to also be even. A first row expansion of $\det Y_n$ with a repeated use the 
induction argument yields:
\begin{align} \nonumber
\det Y_n(p_1,...,p_n) & = (p_1+p_2)\cdot \det Y_{n-1}(p_2,...,p_n) -p_2^2 \cdot \det Y_{n-2}(p_3,...,p_n) \cr 
&= (p_1+p_2)\cdot \sigma_{n-2}(p_2,...,p_n) - p_2^2 \cdot \sigma _{n-3}(p_3,....,p_n) \cr 
& = (p_1+p_2)\cdot ( p_2 \cdot  \sigma_{n-3}(p_3,...,p_n)+ p_3 \cdot ... \cdot p_n)  - p_2^2 \cdot \sigma _{n-3}(p_3,....,p_n) \cr  
& = p_1 \cdot p_2 \cdot \sigma_{n-3}(p_3,...,p_n) + (p_1+p_2) \cdot p_3 \cdot ... \cdot p_n \cr
& = \sigma_{n-1} (p_1,...,p_n) 
\end{align}
completing the proof of the lemma. 
\end{proof}
\section{Computations} \label{examples-section}
In this section we use the results from theorems \ref{main1}, \ref{main2} and \ref{main3} to explicitly evaluate the 
Witt classes of the  knots from examples \ref{example1} -- \ref{example4}. We start with an easy observation. 
\begin{proposition} \label{stab-prop}
If $K$ is a knot obtained from $P(p_1,...,p_n)$ by a finite number of upward stabilizations (see definition \ref{upstab}), then 
$$\varphi (K)= \varphi(P(p_1,...,p_n))$$
Moreover, the signatures of $K$ and $P(p_1,...,p_n)$ are the same and there exists an integer $m$ 
such that $\det K = m^2 \cdot \det P(p_1,...,p_n)$. 
\end{proposition}
This follows easily from theorems \ref{main1}, \ref{main2} and \ref{main3} by inspection. It follows even quicker from observing that the knots $P(p,-p,p_1,...,p_n)$ and $P(p_1,...,p_n)$ are smoothly concordant (see for example \cite{greenejabuka}) and thus in particular also algebraically concordant. This of course implies that their Witt classes are the same and in particular that they have the same signature. Moreover, the determinant  
of a Witt class is well defined up to multiplication by squares.  

We now turn to a more detailed analysis of the examples from section \ref{appandex}. The numerical data presented has already been somewhat simplified by relying on the two relations \eqref{cancelsquare} which 
we use freely and tacitly throughout. 
\vskip3mm
\noindent 
{\bf Example \ref{example1}}
Let $K_1$, $K_2$ and $K_3$ be the knots 
$$K_1=P(21,13,-17,-15,12) \quad \quad K_2=P(-3,-3,-7,5,2) \quad \quad K_3=P(-3,-5,7,9,6)$$
from category $(i)$ in \eqref{categories} and let $K=K_1\#K_2\#K_3$. The $\sigma (K)=0$ but $\varphi (K)$ has order $4$ in $W(\mathbb{Q})$. Thus $K$ has topological and smooth concordance order at least $4$. 

The signatures of $K_1$, $K_2$ and $K_3$ can be computed by a use of theorem \ref{signdet} and are
$\sigma (K_1) = -2$, $\sigma(K_2) =8$ and $\sigma (K_3) = -6$ showing that $\sigma (K) = 0$. The 
rational Witt classes of $K_1$, $K_2$ and $K_3$ are 
\begin{align} \nonumber
\varphi(K_1) & = \langle -56078 \rangle \oplus \langle -105 \rangle \oplus \langle -95 \rangle \oplus 
\langle -38 \rangle \oplus \langle -34 \rangle \oplus  \langle 182 \rangle \oplus \langle 210 \rangle \oplus 
\langle 510510 \rangle \cr
\varphi(K_2) & = \langle 2 \rangle \oplus \langle 2 \rangle \oplus \langle 6 \rangle \oplus 
\langle 6 \rangle \oplus \langle 23 \rangle \oplus  \langle 30 \rangle \oplus \langle 42 \rangle \oplus 
\langle 105 \rangle \cr
\varphi(K_3) & = \langle -42 \rangle \oplus \langle -42 \rangle \oplus \langle -30 \rangle \oplus 
\langle -30 \rangle \oplus \langle -14 \rangle \oplus  \langle -5 \rangle \oplus \langle -3 \rangle \oplus 
\langle -2 \rangle \langle 770 \rangle \oplus \langle 4686 \rangle \cr
\end{align}
Thus, for example, $\partial _{71}(K) = \langle -1 \rangle \in W(\F_{71})\cong \mathbb{Z}_4$ showing that 
$K$ has order $4$ in $W(\mathbb{Q})$.  Similarly, $\partial _{23}(K) = \langle 1 \rangle \in W(\F_{23}) \cong \mathbb{Z}_4$. As a curiosity we note that $\partial _{2549}(K) = \langle 1 \rangle \in W(\F_{2549})\cong 
\mathbb{Z}_2\oplus \mathbb{Z}_2$. 
\vskip3mm
\noindent
{\bf Example \ref{example2}}
Let $K_1$ and $K_2$ be the knots 
$$K_1=P(7,3,-5,2) \quad \quad K_2=P(-19,-15,21,10)$$
from category $(ii)$ in \eqref{categories} and let $K=K_1\#K_2\#K_2\#K_2$. The $\sigma (K)=0$ but $\varphi (K)$ has order $4$ 
in $W(\mathbb{Q})$ and therefore also in the topological and smooth concordance group. 

The signatures of $K_1$ and $K_2$ are found from theorem \ref{signdet} as $\sigma (K_1) = -6$ and $\sigma (K_2) = 2$ and so $\sigma (K) = 0$. The rational Witt classes of $K_1$ and $K_2$ are 
\begin{align} \nonumber
\varphi(K_1) & = \langle -34230\rangle \oplus  \langle -42 \rangle \oplus  \langle -30 \rangle \oplus  \langle -6 \rangle \oplus  \langle -2 \rangle \oplus  \langle -2 \rangle \cr
\varphi(K_2) & = \langle -450870\rangle \oplus  \langle -105 \rangle \oplus  \langle -95 \rangle \oplus  \langle 33 \rangle \oplus  \langle 39 \rangle \oplus  \langle 110 \rangle \oplus  \langle 182 \rangle \oplus  \langle 210 \rangle
\end{align} 
From this one then finds that, for example, $\partial _{3}(K) = \langle 1 \rangle \in W(\F_3) \cong \mathbb{Z}_4$
(as $\partial _3(K_1) = 0$ and $\partial _3 (K_2) = \langle -1 \rangle$). Likewise, $\partial _{163}(K) = 
\langle -1 \rangle \in W(\F_{163})\cong \mathbb{Z}_4$ while $\partial _{113}(K) = \langle a \rangle \in 
W(\F_{113}) \cong \mathbb{Z}_2 \oplus \mathbb{Z}_2$ where $a\in \F_{113}$ is any element which isn't a square. 
\vskip3mm
\noindent
{\bf Example \ref{example3}}
Let $K$ be a knot obtained by a finite number of upward stabilization from either 
$$P(-3,9,15,-5-5) \quad \mbox{ or } \quad P(-3,-5,-11,15,15)$$
from category $(iii)$. 
Then the signature of $K$ is zero, the determinant of $K$ is a square but $\varphi(K) \ne 0 \in W(\mathbb{Q})$. 
Consequently, no such $K$ is slice. 

Note that according to proposition \ref{stab-prop}, it suffices to prove the claims for the two given pretzel knots. From theorem \ref{signdet} we find 
$$ \sigma (K_1) = 0, \, \,  \det K_1 = 75^2 \quad \mbox{ and } \quad  \sigma (K_2) = 0, \, \,  \det K_1 = 135^2$$
and the rational Witt classes of $K_1$ and $K_2$ are 
\begin{align} \nonumber
\varphi(K_1) & = \langle 6 \rangle \oplus  \langle 42 \rangle \oplus  \langle -35 \rangle \oplus  \langle -5 \rangle \cr
\varphi(K_2) & = \langle -2\rangle \oplus  \langle -206 \rangle \oplus  \langle 35535 \rangle \oplus  \langle 345 \rangle 
\end{align} 
This shows that for each $i=1,2$ one obtains $\partial _{3}(K_i)= \langle 1 \rangle \oplus \langle 1 \rangle  \in W(\F_{3})\cong \mathbb{Z}_4$
and similarly $\partial _5(K_i) = \langle 1 \rangle\oplus \langle 2 \rangle \in W(\F_5) \cong \mathbb{Z}_2 \oplus \mathbb{Z}_2$ implying that both knots are non-slice. 
\vskip3mm
\noindent
{\bf Example \ref{example4}}
Let $K_1$, $K_2$ and $K_3$ be the knots  
$$K_1=P(21,13,-17,-15,12) \quad \quad K_2=P(-19,-15,21,10) \quad \quad K_3=P(-15,-7,-7,13,11)$$
from the categories $(i)$, $(ii)$ and $(iii)$ from \eqref{categories} and let $K=K_1\#K_2\#K_3$. Then $\sigma (K) = 0$ but $\varphi (K)$ is of order $4$ in $W(\mathbb{Q})$. 

The signature of $K$ is easily found from theorem \ref{signdet}. The rational Witt classes of $K_1$ and $K_2$ 
have already been computed in examples \ref{example1} and \ref{example2} above while the rational Witt class for $K_3$ is 
\begin{align} \nonumber
\varphi(K_3) & = \langle -22 \rangle \oplus  \langle -5698 \rangle \oplus  \langle 3478 \rangle \oplus  \langle 260474 \rangle 
\end{align} 
From these one arrives at $\partial_{2549} (K) = \langle 1 \rangle \in W(\F_3)\cong \mathbb{Z}_4$, 
$\partial _7(K) = \langle -1 \rangle \in W(\F_7) \cong \mathbb{Z}_4$ and also 
$\partial _{163}(K) = \langle -1 \rangle \in W(\F_{163}) \cong \mathbb{Z}_{4}$. Each of these shows that $K$ has
order 4 in $W(\mathbb{Q})$. 
\section{Proofs of theorems \ref{coro1}, \ref{oddoddeven-theorem}, \ref{coro2} and \ref{main5}} \label{last-section}
This section is devoted to the proofs of theorems listed in the title.  We start with a useful lemma to be used in the subsequent arguments. 
\begin{lemma} \label{auxx1}
Let $\wp$ be a prime number and $p>0$ an odd integer. Write $p = \wp ^\ell \cdot \beta$ 
with $\ell \ge 0$ and $\gcd (\wp, \beta) = 1$. Then 
$$
\partial _\wp (\langle 1\cdot 2 \rangle \oplus \langle 2\cdot 3 \rangle \oplus ... \oplus  \langle (p-1)\cdot p \rangle )= 
\left\{
\begin{array}{cl}
0 & \quad ; \quad \mbox{ if $\ell$ is even} \cr
\langle -\beta \rangle) &  \quad ; \quad \mbox{ if $\ell$ is odd}
\end{array}
\right.$$
where $\partial _\wp : W(\mathbb{Q}) \to W(\F_\wp)$ is the homomorphism between Witt rings from section 
\ref{wittrings}. 
\end{lemma}
\begin{proof}
Assume for the moment that $\wp \ge 3$. 
For any ingeter $2\le k<p$, there are two terms in $\langle 1\cdot 2 \rangle \oplus \langle 2\cdot 3 \rangle \oplus ... \oplus  \langle (p-1)\cdot p \rangle$ containing $k$, namely $\langle (k-1)\cdot k \rangle \oplus \langle k\cdot (k+1) \rangle$. If $\wp$ does not divide $k$, then $\partial _\wp (\langle (k-1)\cdot k \rangle \oplus \langle k\cdot (k+1) \rangle) =0$ by definition of $\partial _\wp$. If $\wp$ does divide $k$, say $k=\wp^\ell \cdot \beta$ with 
$\gcd(\beta, \wp ) =1$, then $k\pm1 \equiv \pm 1 \, (\mbox{mod } \wp)$. Therefore 
$\partial _\wp (\langle (k-1)\cdot k \rangle \oplus \langle k\cdot (k+1) \rangle) =\langle -\beta \rangle \oplus 
\langle \beta \rangle = 0$. The only integer not appearing twice as a factor in this way, is $p$ itself leading to result 
stated by the lemma. 

If $\wp=2$ the result follows in the same manner by pairing up $\langle 1\cdot 2 \rangle \oplus \langle 2\cdot 3 \rangle$, $\langle 3\cdot 4 \rangle \oplus \langle 4\cdot 5 \rangle$ etc. and using the fact that $p$ is odd. 

\end{proof}
\begin{proof}[Proof of theorem \ref{coro1}] Let $K=P(p,q,r)$ be a 3-stranded pretzel knot with
$p,q,r$ odd. Recall from theorem \ref{main3} that the rational Witt class of $K$ is given by 
$$ \varphi(K) = \langle p+q\rangle \oplus \langle (p+q) \det K\rangle$$
where $\det K = pq+pr+qr$. Before proceeding, we first re-write this Witt class in a more symmetric manner
using the relations from theorem \ref{wittpresentation}.  Thus 
\begin{align} \nonumber
\varphi(K) & = \langle p+q\rangle \oplus \langle (p+q)^2r+pq(p+q)\rangle\cr
& =  \langle p+q\rangle \oplus  \langle (p+q)^2r \rangle \oplus  \langle pq(p+q) \rangle \ominus 
 \langle (p+q)^4pqr \det K \rangle \cr
 & = \langle p \rangle \oplus \langle q \rangle \ominus \langle pq(p+q) \rangle \oplus \langle r \rangle \oplus  
 \langle pq(p+q) \rangle \ominus 
 \langle pqr \det K \rangle \cr
 & = \langle p \rangle \oplus \langle q \rangle \oplus \langle r \rangle \ominus 
 \langle pqr \det K \rangle
\end{align}
We shall rely on both of these representations of $\varphi(K)$. 

$\bullet$ Using the first representation for $\varphi(K)$ above, it is easy to see that the rational Witt class 
of $K$ is zero precisely when $\det K = -m^2$ 
for some odd integer $m$. 

$\bullet$ 
Assume now that $|\det K | \equiv 3\,  (\mbox{mod } 4)$ and $\det K<0$. Write 
$\det K = - a_1^{m_1} \cdot ... \cdot a_\ell ^{m_\ell}$ where $a_i$ are positive prime numbers.  
The congruence class of $\det K$ mod $4$ implies that there must be an index $i$ with $a_i  \equiv 3\,  (\mbox{mod } 4)$ and with $m_i$ odd. Write $p,q,r$ as $p=a_i^{\ell_1} \beta _1$, $q=a_i^{\ell_2} \beta _2$ and 
$r=a_i^{\ell_3} \beta _3$ with $\gcd (\beta_{j}, a_i) =1$ and $\ell_j \ge 0$. Similarly, write $\det K = a_i^{m_i} \cdot \beta$ with $\gcd (\beta, a_i) = 1$. According to the parities of $\ell_j$ we have 
$$ \begin{array}{lcl}
(\ell_1,\ell_2,\ell_3) = (\mbox{odd,odd,odd}) & \quad \Longrightarrow \quad & \partial _{a_i} ( \varphi(K) ) = \langle \beta_1\rangle 
\oplus \langle \beta_2\rangle \oplus \langle \beta_3\rangle  \cr
(\ell_1,\ell_2,\ell_3) = (\mbox{odd,odd,even}) & \quad \Longrightarrow \quad & \partial _{a_i} ( \varphi(K) ) = \langle \beta_1\rangle 
\oplus \langle \beta_2\rangle \ominus \langle \beta_1\beta_2\beta_3\beta\rangle  \cr
(\ell_1,\ell_2,\ell_3) = (\mbox{odd,even,even}) & \quad \Longrightarrow \quad &  \partial _{a_i} (\varphi(K)) = \langle \beta_1\rangle 
 \cr
 (\ell_1,\ell_2,\ell_3) = (\mbox{even,even,even}) & \quad \Longrightarrow \quad &  \partial _{a_i} ( \varphi(K) )= \langle  \beta_1\beta_2\beta_3\beta\rangle
\end{array}$$
Since $a_i  \equiv 3\,  (\mbox{mod } 4)$ we know that $W(\F_{a_i})\cong \mathbb{Z}_4$ and so the sum/difference of any 3 generators is again a generator. Thus, in all cases, $ \partial _{a_i} (\varphi(K))$ is a generator of $W(\F_{a_i})$
and is therefore of order 4 (the fact that $\sigma(K) = 0$ follows from the assumption that $\det K<0$). 

$\bullet$ Consider the case of $|\det K| \equiv 1\,  (\mbox{mod } 4)$ and $\det K<0$. Note that every prime $\wp$
congruent to $3\,  (\mbox{mod } 4)$ divides $\det K$ with an even power. Write 
$p=\wp^{\ell_1}\beta_1$, $q=\wp^{\ell_2}\beta_2$ and $r=\wp^{\ell_3}\beta_3$ with $\gcd (\beta_j, \wp) =1$.  
Then for every prime  $\wp \equiv 3\,  (\mbox{mod } 4)$ we obtain  
$$ \begin{array}{lcl}
(\ell_1,\ell_2,\ell_3) = (\mbox{odd,odd,odd}) & \quad \Longrightarrow \quad & \partial _{\wp} ( \varphi(K) ) = \langle \beta_1\rangle 
\oplus \langle \beta_2\rangle \oplus \langle \beta_3\rangle \ominus \langle \beta_1\beta_2\beta_3\beta\rangle  \cr
(\ell_1,\ell_2,\ell_3) = (\mbox{odd,odd,even}) & \quad \Longrightarrow \quad & \partial _{\wp} ( \varphi(K) ) = \langle \beta_1\rangle 
\oplus \langle \beta_2\rangle   \cr
(\ell_1,\ell_2,\ell_3) = (\mbox{odd,even,even}) & \quad \Longrightarrow \quad &  \partial _{\wp} (\varphi(K)) = \langle \beta_1\rangle \ominus  \langle \beta_1\beta_2\beta_3\beta\rangle  \cr
(\ell_1,\ell_2,\ell_3) = (\mbox{even,even,even}) & \quad \Longrightarrow \quad &  \partial _{\wp} ( \varphi(K) )= 0
\end{array}$$
Thus $\partial _{\wp}(\varphi(K))$ is of order $0$ or $2$ in $W(\F_\wp)$. 

$\bullet$ $\varphi(K)$ is of infinite order in $W(\mathbb{Q})$ if and only if $\sigma (K)\ne 0$ which in turn 
occurs if and only if $\det K >0$. 
\end{proof}

The following is a slightly more detailed version of theorem \ref{oddoddeven-theorem}. 
\begin{theorem}  \label{oddoddeven-theorem-general}
Let $K=P(p,q,r)$ with $p,q$ odd and with $r\ne 0$ even. Then 
$\varphi(K)$ is of finite order in $W(\mathbb{Q})$ if and only if 
$$p+q=0\quad \quad \quad \mbox{ or } \quad \quad \quad p+q=\pm2 \, \, \mbox{ and } \, \,\det K>0 $$  
The order of $\varphi(K)$ in $W(\mathbb{Q})$ in these cases is as follows:
\begin{itemize}
\item If $p+q=0$ then $\varphi(K)$ has order $1$ in $W(\mathbb{Q})$. 
\item If $p+q=\pm 2$ and $\det K>0$ then 
$\partial_2(\varphi(K))=0$ and $\partial_\wp(\varphi(K)) = \partial _\wp ( \langle 2 \det K \rangle)$ for every odd prime $\wp$. Consequently
\begin{itemize}
\item  $\varphi(K)$ is of order $1$ in $W(\mathbb{Q})$ if $\det K=m^2$ for some odd integer $m$. 
\item  $\varphi(K)$ is of order $2$ in $W(\mathbb{Q})$ if $\det K$ is not a square and is congruent to
$1 \, (\mbox{mod } 4)$. 
\item  $\varphi(K)$ is of order $4$ in $W(\mathbb{Q})$ if $\det K \equiv 3 \, (\mbox{mod } 4)$.
\end{itemize} 
\end{itemize}
Recall that $\det K =  pq+pr+qr$.
\end{theorem}
\begin{proof} 
The signature of $K$ is zero if and only if (cf. theorem \ref{signdet}):
\begin{itemize}
\item $p+q=0$.
\item $p+q=\pm2$ and $\det K>0$. 
\end{itemize}
In all other cases $\varphi(K)$ is of infinite order in $W(\mathbb{Q})$. If $p+q=0$ then theorem 
\ref{main1} shows that $\varphi(K)=0$ without any condition on $\det K$. 

Turning to the case of $p+q=\pm2$ and $\det K >0$, we first assume, by passing to the mirror image of $K$ if necessary, that $p+q=2$. By interchanging the roles of $p$ and $q$ if needed, we additionally assume that $p>0$. Note that these changes do not affect the sign of $\det K$. The condition $p+q=2$ 
implies that $p>0$ and $q<0$ with the exception of $p=1=q$. We single out this special case first. 
Theorem \ref{main1} shows that the rational Witt class of $K$ in the case of $p=q=1$ is   
$$ \varphi(P(1,1,r)) = \langle -2\rangle \oplus \langle 2\det K\rangle  $$ 
Thus $\partial _2(\varphi(P(1,1,r))) = 0$ and $\partial _\wp(\varphi(P(1,1,r))) = \partial _\wp (\langle2 \det K\rangle )$ for any odd prime $\wp$.

We proceed by keeping our assumptions $p+q=2$, $p>0$ and consider the more general case of $q<0$. Note that 
the rational Witt class of $K$ now takes the form
\begin{align} \nonumber
 \varphi(K) & = \left( \langle -1 \cdot 2\rangle \oplus ... \oplus  \langle  -(p-1)\cdot p \rangle \right)  
 \oplus  \left(  \langle  1 \cdot 2\rangle \oplus ... \oplus  \langle  \cdot (p-3)\cdot (p-2) \rangle \right)\oplus \cr
& \, \oplus  \langle 2p(p-2) \rangle \oplus \langle 2 (\det K) \rangle
\end{align}

Let $\wp$ be a prime number and consider the following cases. 

\begin{enumerate}
\item If $\wp > 2$ and $\wp | (p-2)$, say $p-2= \wp ^\ell \cdot \beta$
with $\gcd (\wp , \beta) = 1$, then with the help of lemma \ref{auxxx1} one obtains 
$$\partial _\wp (\varphi(K)) = \left\{
\begin{array}{cl}
\partial _\wp (\langle 2 \det K)) & \mbox{ if $\ell$ is even } \cr
\langle -\beta \rangle \oplus \langle 2p\beta\rangle \oplus \partial _\wp (\langle 2 \det K))  &  \mbox{ if $\ell$ is odd }
\end{array} \right.
$$
But if $p-2 \equiv 0 \mod \wp$ then $p\equiv 2 \mod \wp$ and so $\langle 2p\beta\rangle = \langle \beta \rangle
\in W(\F_\wp)$. Therefore $\partial _\wp (\varphi(K)) = \partial _\wp (\langle 2 \det K)) $. 
\item If $\wp > 2$ and $\wp | p$, say $p= \wp ^\ell \cdot \beta$ with $\gcd (\wp, \beta)=1$, then using again lemma 
\ref{auxxx1} we get 
$$\partial _\wp (\varphi(K)) = \left\{
\begin{array}{cl}
\partial _\wp (\langle 2 \det K))  & \mbox{ if $\ell$ is even } \cr
\langle \beta \rangle \oplus \langle 2\beta (p-2)\rangle \oplus \partial _\wp (\langle 2 \det K)) &  \mbox{ if $\ell$ is odd }
\end{array} \right.
$$
But  $p \equiv 0 \mod \wp $ implies that $p-2\equiv -2 \mod \wp $ so that 
$\langle 2\beta(p-2)\rangle = \langle -4\beta \rangle = \langle -\beta \rangle \in W(\F_\wp)$. Thus we get 
again that $\partial_\wp (\varphi(K)) = \partial _\wp (\langle 2 \det K)) $. 
\item If $\wp > 2$ and $\wp$ doesn't divide either of $p$ or $p-2$,  then $\partial _\wp (\varphi(K))$ is trivially equal to $\partial _\wp (\langle 2 \det K)) $ (with the help of lemma \ref{auxxx1}). 
\item Consider $\wp =2$. Since $p$ and $p-2$ are odd, it is easy to see that 
the determinant $\det \varphi(K)$ is of the form $2^{2\ell} \cdot \beta$ for some odd $\beta$. But then 
$\partial _2(\varphi(K)) = 0$ by definition.  
\end{enumerate}
Thus we obtain $\partial _\wp (\varphi(K)) = \partial _\wp (\langle 2 \det K))$ for all odd prime integers $\wp$
and $\partial _2 (\varphi(K)) = 0$. 

Given this, it is now an easy matter to verify the stated orders of $\varphi(K)$ in $W(\mathbb{Q})$. For example, 
if $\det K=m^2$ then $\partial _\wp(\langle 2\det K\rangle)=0$ for all primes $\wp$ and thus $\varphi(K)=0 \in W(\mathbb{Q})$. If $\det K \equiv 3 \, (\mbox{mod } 4)$ then there must exist a prime $\wp  \equiv 3 \, (\mbox{mod } 4)$ dividing $\det K$ with an odd power. Therefore $\partial _\wp (\langle 2 \det K\rangle)$ 
yields a generator of $W(\mathbb{\F_\wp})\cong \mathbb{Z}_4$. We leave the remaining case as an easy exercise for the interested reader. 
\end{proof}
In preparation for the proof of theorem \ref{coro2}, we state a couple of auxiliary lemmas first. 
\begin{lemma} \label{halo1}
Consider odd integers $n,p_1,...,p_{n-1}$ with $n\ge 3$ and let $p_n\ne 0$ be an even integer. 
Let $\wp$ be an odd prime which doesn't divide any of $p_1,...,p_{n-1}$ and  assume that 
$\det P(p_1,...,p_n) = \pm m^2$ for some integer $m$. Then $\partial _\wp (\varphi(P(p_1,...,p_n))) = 0$.
\end{lemma}
\begin{proof}
There are two cases which we consider separately, namely the case when $\wp$ divides $\det P(p_1,...,p_n)$
and the case when it doesn't. Let us write $\det P(p_1,...,p_{n}) = \varepsilon \cdot m^2$ for some 
choice of $\varepsilon \in \{\pm 1\}$. 

Assume firstly that $\wp$ is a divisor of $\det P(p_1,...,p_n)$. 
By lemma \ref{auxx1} and theorem \ref{main1} we find that 
\begin{align} \nonumber
 \partial _\wp (\varphi(P(p_1,...,p_n))) & = \partial _\wp ( \left\langle -(p_1\cdot ... \cdot p_{n-1}) \cdot \det P(p_1,...,p_{n-1}) \right\rangle) \oplus \cr
 &\quad \quad \quad \quad \quad \quad \quad \quad \quad\oplus \partial _\wp (\left\langle \det P(p_1,...,p_{n-1}) \cdot \varepsilon \right\rangle) 
 \end{align}
Since 
\begin{equation} \label{detrelations}
\det P(p_1,...,p_n) = p_n \cdot \det P(p_1,...,p_{n-1}) +  p_1\cdot ... \cdot p_{n-1}
\end{equation}
and $\wp$ divides $\det P(p_1,...,p_n)$ but does not divide $p_1\cdot ... \cdot p_{n-1}$, we see that $\wp$ cannot divide $p_n \cdot \det P(p_1,...,p_{n-1})$. Thus $\partial _\wp (\varphi(P(p_1,...,p_n)))   = 0$. 

Next, suppose that $\wp$ does not divide $\det P(p_1,...,p_n)$. Write $\det P(p_1,...,p_{n-1}) = \wp ^\ell \cdot \beta$ for some integer $\ell \ge 0$ and some $\beta$ with $\gcd ( \wp, \beta) = 1$. If $\ell$ is even then 
$\partial _\wp (\varphi(P(p_1,...,p_n)))$ vanishes trivially. Else, if $\ell$ is odd, and using \eqref{detrelations} again, 
we see that $\varepsilon \cdot p_1\cdot ... \cdot p_{n-1}$ is a square modulo $\wp$. Therefore, 
\begin{align} \nonumber
\partial _\wp (\varphi(P(p_1,...,p_n))) & = \langle - (p_1\cdot ... \cdot p_{n-1}) \cdot \beta \rangle \oplus \langle 
\varepsilon \cdot  \beta \rangle  \cr
& = \langle -\varepsilon \cdot \beta \rangle \oplus \langle \varepsilon \cdot  \beta \rangle  \cr
& = 0  
\end{align}
\end{proof}
\begin{lemma}  \label{halo2}
Consider again odd integers $n,p_1,...,p_{n-1}$ with $n\ge 3$ and let $p_n\ne 0$ be an even integer. 
Let $\wp$ be an odd prime which divides exaclty one $p_i \in \{p_1,...,p_{n-1}\}$. 
Assume again that $\det P(p_1,...,p_n) = \pm m^2$ for some integer $m$. Then $\partial _\wp (\varphi(P(p_1,...,p_n))) = 0$.
\end{lemma}
\begin{proof}
For concreteness assume that $\wp$ divides $p_1$ and that therefore $\gcd(\wp,p_j) =1$ for all $j=2,...,n-1$. 
The assumption $\det P(p_1,...,p_n) = \pm m^2$ along with lemma \ref{auxx1} and theorem \ref{main1}, implies that  
\begin{align} \nonumber
 \partial _\wp (\varphi(P(p_1,...,p_n))) & = \partial _\wp (\langle Sign(-p_1) (|p_1|-1)\cdot  |p_1| \rangle) \oplus  \cr
 & \quad \quad \quad \quad  \oplus  \partial _\wp (\left\langle -(p_1\cdot ... \cdot p_{n-1}) \cdot \det P(p_1,...,p_{n-1}) \right\rangle) \oplus \cr 
 & \quad \quad \quad \quad \oplus \partial _\wp (\left\langle \pm \det P(p_1,...,p_{n-1}) \right\rangle)
 \end{align}
Since 
$$\det P(p_1,...,p_{n-1}) =  p_1\cdot \left( \prod _{i=2}^{n-1} p_2\cdot ... \cdot \hat p_i \cdot ... \cdot p_{n-1}\right) + 
p_2\cdot ... \cdot p_{n-1} $$
we see that $\wp$ cannot divide $\det P(p_1,...,p_{n-1})$, in fact, 
$$\det P(p_1,...,p_{n-1}) \equiv  p_2\cdot ... \cdot p_{n-1} \mbox{ (mod $\wp$)}$$ 
Let us write $p_1 = \wp ^\ell \cdot \beta$ for some $\ell \ge 0$ and with $\gcd (\wp , \beta) = 1$. If $\ell $ is odd, then  
\begin{align} \nonumber
 \partial _\wp (\varphi(P(p_1,...,p_n))) & = \langle \beta \rangle \oplus \left\langle -\beta \cdot p_2 \cdot  ... \cdot p_{n-1} \cdot \det P(p_1,...,p_{n-1}) \right\rangle \cr 
 & = \langle \beta \rangle \oplus \left\langle -\beta \cdot (p_2 \cdot  ... \cdot p_{n-1})^2  \right\rangle \cr
 & = \langle \beta \rangle \oplus \left\langle - \beta \right\rangle \cr
 & = 0 
\end{align}
On the other hand, if $\ell $ is even, then $ \partial _\wp (\varphi(P(p_1,...,p_n)))=0$ by the definition of the map 
$ \partial _\wp$. 
\end{proof}
The results from lemmas \ref{halo1} and \ref{halo2} imply the statement of theorem \ref{coro2}. 
\begin{proof}[Proof of theorem \ref{main5}]
We start by finding the linking matrix $\mathcal L$ of $K=P(5,-3,8)$ as in section \ref{one-one}. The formulae 
provided there easily imply that 
$$\mathcal L = \left[
\begin{array}{rrrr|rr|rr}
-1 & -1& -1& -1&0 &0 &-1 &0 \cr
0 & -1& -1& -1&0 &0 &-1 &0 \cr
0 & 0& -1& -1&0 &0 &-1 &0 \cr
0 & 0& 0& -1&0 &0 &-1 &0 \cr \hline 
0 & 0& 0& 0&1 &1 & 1&0 \cr
0 & 0& 0& 0&0 &1 & 1&0 \cr \hline 
 0& 0&0 & 0& 0 &0 &0  & 0\cr
0 & 0& 0& 0&0 &0 & 1 & 4 
\end{array}
\right]
$$
Pick $\omega =a+ib \in S^1\subset \mathbb{C}$ (so that $a^2+b^2=1$) and form the matrix $A_\omega = (1-\omega)\mathcal L + (1-\bar\omega)\mathcal L^\tau$. By definition, the Tristram-Levine signature $\sigma _\omega (K)$ of $K$ equals the signature of $A_\omega$. It is well known that the signatures $\sigma _\omega (K)$ are 
constant away from the unit roots of the symmetric Alexander polynomial $\Delta _K(t)$. We thus turn to computing the latter.  

The Alexander polynomial $\Delta_K(t) = \det ({t}^{1/2} \mathcal L - t^{-1/2} \mathcal L^\tau )$ of $K=P(5,-3,8)$  is given by $\Delta_K(t) = t^3-2t^2-t+5-t^{-1}-2t^{-2}+t^{-3}$. Its graph is depicted in figure \ref{pic8}. 
\begin{figure}[htb!] 
\centering
\includegraphics[width=10cm]{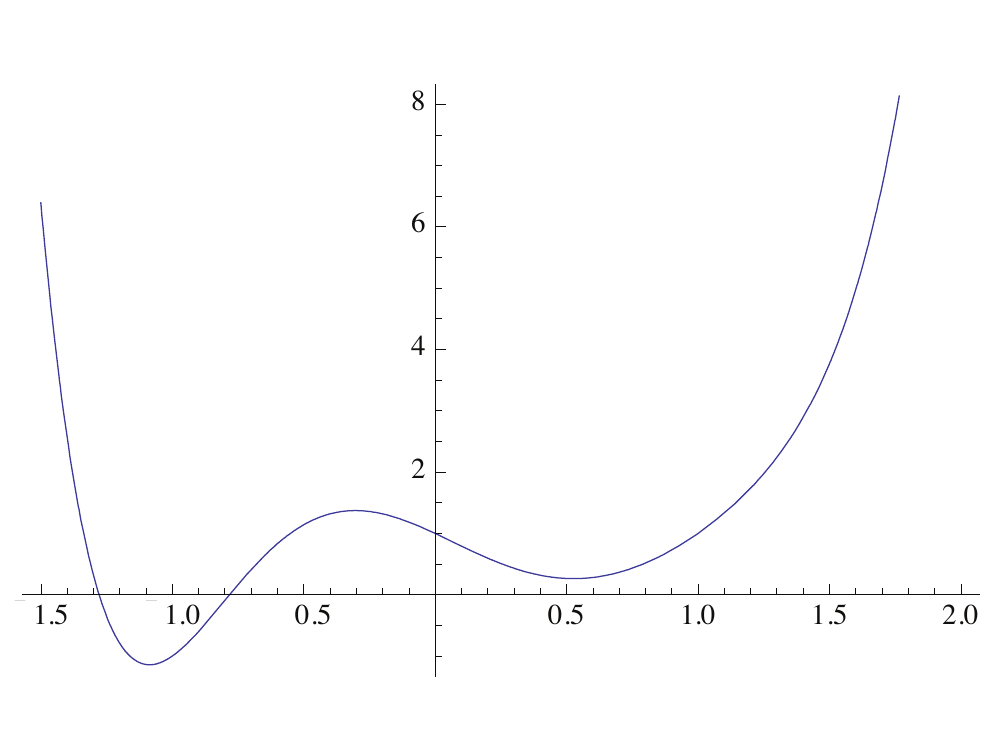} 
\caption{The graph of $t^3\cdot \Delta_{P(5,-3,8)}(t)$.  }  \label{pic8}
\end{figure}
Clearly visible on the graph, the two real roots $t_{1,2}$ of $\Delta _K(t)$ are not of unit norm. The 4 complex roots 
are approximately 
$$ t_{3,4} = 0.528853 \pm 0.269329\,  i \quad \quad \mbox{ and } \quad \quad t_{5,6} = 1.50147\pm 0.764653 \, i $$
showing that the approximate norms of $t_{3,4}$ and $t_{5,6}$ are 
$$ |t_{3,4}| = 0.352223 \quad \quad \mbox{ and } \quad \quad  |t_{5,6}| = 2.83911$$
Thus $\Delta_K(t)$ has no roots on $S^1$ so that $\sigma_\omega (K) = \sigma(K)$ for all $\omega \in S^1$. 
But $\sigma (K)=0$ as is easily computed from theorem \ref{signdet}. This implies that $K$ is of finite algebraic concordance order, cf. \cite{chuck}. 

On the other hand, if $K$ were algebraically slice, then we could factor $\Delta _K(t)$ as $f(t) \cdot f(t^{-1})$ 
for some $f(t) \in \mathbb{Z}[t]$. This however is not the case. An easy way to see this is to note that the 
mod 2 reduction of $\Delta_K(t)$ looks like 
\begin{align} \nonumber
\Delta_K(t) & \equiv  t^3+t+1+t^{-1}+t^{-3} (\mbox{ mod } 2) \cr
& \equiv  (t+1+t^{-1})(t^2+t+1+t^{-1}+t^{-2})   (\mbox{ mod } 2) \cr
\end{align}
Now, $t+1+t^{-1}$ is irreducible in $\mathbb{Z}_2[t,t^{-1}]$ but $t^2+t+1+t^{-1}+t^{-2}$ is not divisible by $t+1+t^{-1}$. Thus $\Delta_K(t)$ could not have factored as $f(t) \cdot f(t^{-1})$ and so $K$ is not algebraically slice. 
In fact, using {\sc Mathematica} one finds that $\Delta_K(t)$ is  irreducible over $\mathbb{Z}[t,t^{-1}]$. 

Finally, the fact that $\varphi(K)=0\in W(\mathbb{Q})$ follows readily from theorem \ref{oddoddeven-theorem} since $\det K = 1$, $5+(-3)=2$
and, as already mentioned, $\sigma (K)=0$. 
\end{proof}

\end{document}